\pdfoutput=1
\RequirePackage{ifpdf}
\ifpdf 
\documentclass[pdftex]{sigma}
\else
\documentclass{sigma}
\fi

\usepackage{tikz}
\usepackage[all]{xy}

\usetikzlibrary{arrows.meta}

\usetikzlibrary{shapes.misc}
\usetikzlibrary{decorations.pathmorphing}

\usepackage{euscript}

\DeclareMathOperator{\Aut}{Aut}
\DeclareMathOperator{\Bir}{Bir}
\DeclareMathOperator{\Auteq}{Auteq}
\DeclareMathOperator{\Ext}{Ext}
\DeclareMathOperator{\Quad}{Quad}

\newcommand {\Coh} {\operatorname{Coh}}
\newcommand {\ccoh} {\operatorname{coh}}

\newcommand{\Spec}{\mathrm{Spec}}

\newcommand{\Ham}{\mathrm{Ham}}
\newcommand{\FS}{\mathrm{FS}}
\newcommand{\std}{\mathrm{std}}
\newcommand{\Bl}{\mathrm{Bl}}
\newcommand{\Adm}{\mathrm{Adm}}

\def\scrA{\EuScript{A}}
\def\scrJ{\EuScript{J}}

\def\scrQ{\EuScript{Q}}
\def\scrR{\EuScript{R}}

\def\cL{\EuScript{L}}
\def\cO{\EuScript{O}}
\def\cT{\EuScript{T}}
\newcommand{\scrF}{\EuScript{F}}
\newcommand{\scrC}{\EuScript{C}}
\newcommand{\scrE}{\EuScript{E}}
\def\scrW{\EuScript{W}}

\DeclareMathOperator{\Diff}{Diff}

\DeclareMathOperator{\Symp}{Symp}
\DeclareMathOperator{\Stab}{Stab}

\def\bK{\mathbb{K}}
\def\bP{\mathbb{P}}

\def\bG{\mathbb{G}}
\def\bZ{\mathbb{Z}}
\def\bQ{\mathbb{Q}}

\def\bC{\mathbb{C}}

\def\Db{\mathop{\rm{D}^b}\nolimits}

\def\bT{\mathbb{T}}
\def\bR{\mathbb{R}}

\def\A{A}
\def\scrD{\EuScript{D}}

\def\scrS{\EuScript{S}}

\DeclareMathOperator{\Pic}{Pic}

\newcommand{\id}{\mathrm{id}}

\usepackage{mathabx,graphicx}

\def\Circlearrowright{\ensuremath{%
 \rotatebox[origin=c]{180}{$\circlearrowright$}}}

\newtheorem{thm}{Theorem}[section]
\newtheorem{lem}[thm]{Lemma}
\newtheorem{prop}[thm]{Proposition}
\newtheorem{cor}[thm]{Corollary}
\newtheorem{qn}[thm]{Question}
\newtheorem{conj}[thm]{Conjecture}
\newtheorem{assum}[thm]{Assumption}

\theoremstyle{definition}

\newtheorem{rmk}[thm]{Remark}

\DeclareMathOperator{\tw}{\mathsf{tw}}

\newcommand{\EuI}{\EuScript{I}}

\numberwithin{equation}{section}

\begin{document}

\allowdisplaybreaks

\newcommand{\arXivNumber}{2601.20438}

\renewcommand{\thefootnote}{}

\renewcommand{\PaperNumber}{080}

\FirstPageHeading

\ShortArticleName{Splitting Symplectic Monodromy}

\ArticleName{Splitting Symplectic Monodromy\footnote{This paper is a~contribution to the Special Issue on Geometry and Dynamics in memory of Will Merry. The~full collection is available at \href{https://sigma-journal.com/Merry.html}{https://sigma-journal.com/Merry.html}}}

\Author{Ailsa KEATING~$^{\rm a}$, Ivan SMITH~$^{\rm a}$ and Michael WEMYSS~$^{\rm b}$}
\AuthorNameForHeading{A.~Keating, I.~Smith and M.~Wemyss}

\Address{$^{\rm a)}$~Centre for Mathematical Sciences, University of Cambridge,\\
\hphantom{$^{\rm a)}$}~Wilberforce Road, CB3 0WB, UK}
\EmailD{\mail{amk50@cam.ac.uk}, \mail{is200@cam.ac.uk}}
\Address{$^{\rm b)}$~School of Mathematics and Statistics, University of Glasgow,\\
\hphantom{$^{\rm b)}$}~University Place, Glasgow G12 8QQ, UK}
\EmailD{\mail{michael.wemyss@glasgow.ac.uk}}

\ArticleDates{Received January 29, 2026, in final form August 06, 2026; Published online August 24, 2026}

\Abstract{We discuss some examples in which symplectic monodromy (provably or conjecturally) splits off the symplectic mapping class group, hoping to illustrate different techniques and inputs to the arguments. Along the way we formulate several open questions and conjectures. \emph{Inter alia}, we also construct a compact symplectic six-manifold which contains infinitely many pairwise disjoint Lagrangian 3-spheres.}


\Classification{53D37; 14J33; 53D99}

\renewcommand{\thefootnote}{\arabic{footnote}}
\setcounter{footnote}{0}

\section{Introduction}\label{sec:intro}

Let $(X,\omega)$ be a symplectic manifold; we will typically be interested in non-compact manifolds well-behaved (e.g.,~convex) at infinity. Suppose $X$ moves in a moduli space $\mathcal{M}_X$ which carries a universal family $\scrC \to \mathcal{M}_X$ which is a locally Hamiltonian fibre bundle. Symplectic parallel transport gives rise to a monodromy representation
\begin{equation*} 
\rho\colon \ \pi_1(\mathcal{M}_X) \to \pi_0\Symp(X,\omega).
\end{equation*}
This may be valued in $\pi_0\Symp_{\rm ct}(X,\omega)$ depending on the geometry of the family at infinity. Such representations provide arguably the primary source of constructions of interesting symplectomorphisms. This suggests that a fundamental question is to understand the map $\rho$.

To understand a group, one should make it act on something. The group $\pi_0\Symp(X,\omega)$ acts by autoequivalences of the (derived) Fukaya category $\scrF(X)$, defining a representation
\begin{equation*} 
\pi_0\Symp(X,\omega) \to \mathrm{Auteq}(\scrF(X)).
\end{equation*}
This may be up to shift; we will consider gradings, and the question of compact or non-compact support of maps, more carefully in the main text.
\emph{A priori} the symplectic topology of $X$ has no reason to know anything about its incarnation as a special fibre of $\scrC \to \mathcal{M}_X$ (for instance, its incarnation as an algebro-geometric complete intersection), and yet remarkably, in a few non-trivial cases, one can relate the autoequivalence group to $\pi_1(\mathcal{M}_X)$ itself, and infer that $\rho$ \emph{splits}.

This note illustrates that phenomenon in some basic cases: certain Milnor fibres (Theorem~\ref{thm:milnor case}), quiver 3-folds (Theorem~\ref{thm:bordered surface case}) and mirrors to 3-fold flopping contractions (Theorem~\ref{thm:double bubble case}).
The results in Section~\ref{sec:birational-proj-plane}
have a different character, morally arising from monodromy in a space of integral affine structures on an SYZ base.
 The point of gathering examples together is to emphasise that it is an interesting question to understand how widespread this phenomenon is, and to show that establishing it in any given setting typically requires non-trivial machinery. The extent to which the collection of cases gathered here is prototypical is itself an interesting question.

\begin{rmk}
The paper does not discuss an exhaustive list of splittings of symplectic monodromy, but we gather a few other examples in Section~\ref{sec:sporadic}, for completeness and to help establish the broader scene for the reader.
\end{rmk}

Besides the main `splitting' theme, the paper contains a few digressions. One we draw attention to here is in Section~\ref{sec:conifold}, where we construct a compact symplectic six-manifold which contains infinitely many pairwise-disjoint Lagrangian 3-spheres. The existence of such examples was a well-known open question in the early days of the subject.

\bigskip

\textbf{Outline of the paper.} To aid the reader:
\begin{enumerate}\itemsep=0pt
\item Section~\ref{Sec:flavours} recalls the compact, wrapped and Rabinowitz Fukaya categories.
 \item Sections \ref{sec:type-A-Milnor-fibres}, \ref{sec:marked-bordered-surfaces} and~\ref{sec:double-bubble} unconditionally prove splitting results in three cases, respectively: $A_k$\nobreakdash-Milnor fibres of any dimension (i.e., linear plumbings of cotangent bundles of spheres at points); non-compact Calabi--Yau 3-folds associated to marked bordered surfaces (which generalise 3-dimensional $A_k$-Milnor fibres, corresponding to the bordered surface being a~disc); and `multibubble' plumbings (i.e.,~linear plumbings of 3-spheres along circles rather than along points). Each of these cases relies on a particular classification result: for the exchange graph on the corresponding cluster (in this case Rabinowitz) category (Section~\ref{sec:type-A-Milnor-fibres}); for the class of Lagrangian spheres in the 3-fold (Section~\ref{sec:marked-bordered-surfaces}); and for all bounded hearts in the derived category of a suitable algebro-geometric model for the compact Fukaya category (Section~\ref{sec:double-bubble}).
 \item Section~\ref{sec:conifold} recalls a splitting in the case of the `conifold mirror' (a double bubble cousin given as a certain plumbing of two $T^*S^3$s along a Hopf link) and then discusses inexact deformations of this space, and formulates a question about splitting in that context. The obstruction to answering that is essentially the lack of a classification of Lagrangian spheres in the inexact setting.
 \item Section~\ref{sec:birational-proj-plane} discusses the birational mirror to the projective plane, and proves a splitting result conditional on a certain mirror symmetry statement (extending the known case of exact HMS to a non-exact setting). This conjecture is probably within reach of known technology, but would be a major digression, so we axiomatise our requirements.
 \item Section~\ref{sec:sporadic} gathers various additional results and questions, including a conjectural \emph{non-splitting} result for certain projective complete intersections. It seems new ideas would be needed to prove such a statement.
\end{enumerate}

\textbf{Dedication.}
 Will Merry spent several years at Cambridge, bringing his energy and warmth and enthusiasm for mathematics to the symplectic geometry group and beyond. One of us (I.S.) was his PhD thesis examiner, whilst another (A.K.) first encountered symplectic manifolds through him. We continue to benefit both from his research and from his mentorship of younger mathematicians throughout his time at ETH. His work includes fundamental contributions to Rabinowitz Floer theory and the behaviour of Weinstein manifolds under inexact deformations, each of which make appearances in this paper. We are honoured to contribute this paper to a~volume dedicated to his memory.

\section{Flavours of Fukaya category}\label{Sec:flavours}

Let $(X,\omega)$ be a connected symplectic manifold with $2c_1(X)=0$. We will fix a Maslov cover of the Lagrangian Grassmannian, meaning that there is a well-defined group $\Symp^{\rm gr}(X)$ of graded symplectomorphisms, and a subgroup $\bZ \leq \Symp^{\rm gr}(X)$ generated by the shift (i.e., the fibrewise deck transformation). Recall from \cite[Lemma~2.4]{Seidel:graded} the following.

\begin{lem}
There is an exact sequence
\[
1 \to \bZ \to \Symp^{\rm gr}(X) \to \Symp(X) \to H^1(X;\bR).
\]
The final map is not a homomorphism but its kernel is a subgroup.
\end{lem}

\begin{lem}\label{lem:forgetting-gradings-splits}
Suppose that $X$ is non-compact. The forgetful map $\Symp^{\rm gr}(X) \to \Symp(X)$ is split over the subgroup $\Symp_{\rm ct}(X)$.
\end{lem}

\begin{proof} Take the unique gradings for symplectomorphisms which act trivially on the fibres of the Maslov cover near infinity.
\end{proof}

Let $(X,\omega)$ be a Liouville manifold, with contact type boundary, and assume $2c_1(X)=0$. We will be interested in all of
\begin{enumerate}\itemsep=0pt
\item[(1)] the compact Fukaya category $\scrF(X)$, with objects compact graded exact spin Lagrangians;
\item[(2)] the wrapped Fukaya category $\scrW(X)$, with objects graded exact spin Lagrangians which are cones over Legendrian submanifolds near infinity;
\item[(3)] the Rabinowitz Fukaya category $\scrR(X) := \scrW(X) / \scrW(X)^{\rm prop}$, which is the Verdier quotient of the (Yoneda image of the) wrapped category by its subcategory of proper modules.
\end{enumerate}

All three categories are defined over $\bZ$, and $\scrF(X) \subset \scrW(X)$ is a full and faithful subcategory. We say that `Koszul duality' holds for the compact and wrapped categories if $\scrF(X) \subset \scrW(X)$ is a Koszul dualising subcategory in the sense of \cite[Definition~2.12]{GGV}. This often holds for Weinstein manifolds obtained from plumbings of simply connected pieces. Concretely, the wrapped category of a Weinstein manifold admits a finite set of generators $\{T_i \mid 1 \leq i \leq N\}$ and one seeks a finite set of compact branes $\{L_i \mid 1 \leq i \leq N\}$ such that the total endomorphism algebras $\scrQ = \bigoplus_{i,j} HF^*(L_i,L_j)$ and $\scrA = \bigoplus_{i,j} HW^*(T_i,T_j)$, which are augmented over the semisimple ring $\Bbbk = \bigoplus_{i=1}^N \bK_i$ generated by the units of the objects, are $A_{\infty}$-Koszul dual, meaning
\[
\scrQ \simeq \Ext_{\scrA}(\Bbbk,\Bbbk), \qquad \scrA \simeq \Ext_{\scrQ}(\Bbbk,\Bbbk).
\]
In this setting, there is an equivalence \cite[Corollary~1.4]{GGV}
\[
\scrR(X) \simeq \scrW(X) / \scrF(X).
\]
This is the formulation we will encounter in practise, but the general algebraic construction of~$\scrR(X)$ from~$\scrW(X)$ in particular shows that autoequivalences of the latter act on the former.
 (In Section~\ref{sec:birational-proj-plane}, it will be important to work over~$\bC$, to incorporate $\bC^*$-valued local systems on Lagrangians.)

 Any autoequivalence of $\scrW(X)$ arising from a symplectomorphism of $X$ will preserve the compact category $\scrF(X)$, and hence induce an equivalence of $\scrR(X)$. For background on $\scrF(X)$, see \cite{Seidel:FCPLT}, and for background on $\scrW(X)$, see~\cite{Abouzaid-Seidel:Viterbo}; for background on derived localisation in general, see~\cite{Drinfeld}, and in the context of $\scrR(X)$ specifically, see~\cite{GGV}.

The Grothendieck group $K(\scrF(X))$ is the free abelian group on isomorphism classes of objects, modulo relations coming from exact triangles.
For later use, we define the `numerical Grothendieck group'
\begin{equation}\label{eqn:Knum}
K_{\rm num}(\scrF(X))
\end{equation}
to be the quotient of $K(\scrF(X))$ by the subgroup which is the radical of the Euler pairing, i.e., that generated by compact exact Lagrangians $L \in \mathrm{Ob}\,\scrF(X)$ for which $\chi(HF^*(L,\bullet)) = 0$ for every $\bullet \in \scrW(X)$.

We will be interested in actions of symplectomorphism groups on Fukaya categories. Especially for actions on the wrapped category, one needs some control of maps at infinity to ensure that (perhaps after isotopy) they preserve the class of Lagrangians which are conical\footnote{Or more generally geometrically bounded for a dissipative almost complex structure in the sense of Groman~\cite{Groman}.} at infinity. We say a symplectomorphism $\phi$ of a Liouville manifold $(X,{\rm d}\theta)$ is \emph{exact} if $\phi^*\theta = \theta + {\rm d}f$ for some $f\colon X \to \bR$, and \emph{strongly exact} if one can take $f$ to have compact support. Both notions have obvious graded analogues. It is classical that strongly exact maps act on the wrapped category. Write $\Symp_{\rm ex}(X)$ for the group of exact symplectomorphisms.

\begin{lem}\label{lem:acts}
 Assume that the Liouville manifold $X$ has finite type, so can be written as the completion of a Liouville domain $\mathring{X}$. Then
 \begin{enumerate}\itemsep=0pt
 \item[$(1)$] There is a homomorphism \smash{$\pi_0\Symp_{\rm ct}(X) \to \pi_0\Symp_{\rm ex}\bigl(\mathring{X}\bigr)$};
 \item[$(2)$] The group \smash{$\pi_0\Symp\bigl(\mathring{X}\bigr)$} acts on both $\scrF(X)$ and $\scrW(X)$.
 \end{enumerate}
\end{lem}

\begin{proof}
 For the first statement, conjugation by a suitably large Liouville flow will bring the support into $\mathring{X}$. We claim we can then deform to obtain an exact representative of the map. Let $\omega = {\rm d} \lambda$ be the Liouville form on $X$, and suppose we are given $\psi \in \Symp_{\rm ct} (X)$. We have $\psi^\ast \lambda - \lambda = \theta$, where $\theta$ is a closed one-form with compact support. Let~$V$ be the vector field which is symplectically dual to $-\theta$. This has compact support, and so is integrable for all time; let $\sigma_t$ denote its flow. Then applying Cartan's formula, we get that $\psi \circ \sigma_1 $ is an exact symplectomorphism. (Note that the proof of \cite[Lemma~1.1]{BEE} specialises to this argument in this case.)

 This construction can be performed in families, which shows that the inclusion
 \[
 \Symp_{\rm ct,ex}(X) \hookrightarrow \Symp_{\rm ct}(X)
 \]
 is a homotopy equivalence, so induces an isomorphism on $\pi_0$. On the other hand, there is then a forgetful map
 \[
 \Symp_{\rm ct,ex}(X) \to \Symp_{\rm ex}(\mathring{X})
 \]
 which induces a homomorphism on $\pi_0$. Combining these two gives the first statement.

 The second statement is exactly \cite[Corollary~2.4]{Keating-Smith} (cf.~\cite{KS:corr}, which in particular involves deforming a map of the Liouville domain to be exact, using a version of the argument given in the previous paragraphs).
\end{proof}

\begin{lem}
In the setting of Lemma~{\rm \ref{lem:acts}}, assume moreover that $2c_1(\mathring{X})=0$.
Then there are compatible representations
\[
\xymatrix{
\pi_0\Symp^{\rm gr}_{\rm ct}(X) \ar[rr] \ar[d] && \Auteq(\scrF(X)) \\
\pi_0\Symp^{\rm gr}(\mathring{X}) \ar[rr] && \Auteq(\scrW(X)). \ar[u]
}
\]
\end{lem}
The top line is classical; the bottom line is \cite[Corollary~2.9]{Keating-Smith} (see also~\cite{KS:corr}).

Although the Rabinowitz category is not a priori determined by the contact boundary of $X$, one still has the following lemma.

\begin{lem}\label{lem:trivial on Rabinowitz}
The group $\pi_0\Symp^{\rm gr}_{\rm ct}(X)$ acts trivially on the Rabinowitz category.
\end{lem}

\begin{proof}
Let $L$, $K$ be Lagrangians which are cones over Legendrian submanifolds $L_{\infty}$ and $K_{\infty}$.
The morphism group ${RFC}^*(L,K)$ in $\scrR(X)$ between $L$ and $K$ is the chain-level mapping cone of a natural continuation map ${CW}_*(L,K) \to {CW}^*(L,K)$ coming from a monotone homotopy of Hamiltonians from `infinite negative wrapping to infinite positive wrapping'. This factors as a composition
\[
\xymatrix{
{CW}_*(L,K) \ar[r]_-a & {CF}^*(L,\phi_{-\varepsilon}K) \ar[r]_-c &
{CF}^*(L,\phi_{+\varepsilon}K) \ar[r]_-b & {CW}^*(L,K),
}
\]
where $a$, $b$ are projection respectively inclusion maps of summands into the relevant direct systems and $c$ is the natural continuation from small negative to small positive wrapping. Viewing ${RFC}^*(L,K) = \mathrm{Cone}(b\circ c \circ a)$,
there are exact triangles of mapping cones
\[
\mathrm{Cone}(c\circ a) \to {RFC}^*(L,K) \to \mathrm{Cone}(b)
\]
and
\[
\mathrm{Cone}(a) \to \mathrm{Cone}(c\circ a) \to \mathrm{Cone}(c).
\]
The continuation map $c$ induces the Poincar\'e duality isomorphism in Lagrangian Floer theory
\[
HF^*(L,\phi_{-\varepsilon}K) \to HF^*(L,\phi_{+\varepsilon}(K))
\]
(which recovers the isomorphism $H_*(L,\partial_{\infty}L) \to H^{n-*}(L)$ when $L=K$),
so has acyclic cone. On the other hand, $\mathrm{Cone}(b) \simeq {CW}(L,K)^+$ defines the Lagrangian Floer analogue of `positive symplectic cohomology' \cite{Ritter}; it can be computed by an action truncated version of the wrapped Floer complex which removes action zero chords (i.e., interior intersection points), yielding a~complex only built from Reeb chords at infinity between the Legendrian boundaries. Similarly for $\mathrm{Cone}(a)$. Putting these ingredients together, one obtains a cochain complex quasi-isomorphic to that underlying ${RFC}^*(L,K)$ and which is generated only by Reeb chords at the boundary. Any compactly supported map then acts trivially on the underlying chain group.

 Since $L$ and $K$ were arbitrary, compactly supported maps act trivially on the whole $A_{\infty}$-category.\looseness=1
\end{proof}

To conclude this section, we briefly recall the categories associated to quivers with potential; a more extensive discussion can be found in~\cite{Smith:quiver}.

A quiver $Q$ is a directed graph, and a potential $W$ on $Q$ is a formal sum of cycles in the edges. To such $(Q,W)$ one can associate a 3-dimensional\footnote{When the arrows admit suitable gradings, one can associate a ${\rm CY}_n$-category for any $n$.} Calabi--Yau category $\scrC(Q,W)$ via the `Ginzburg construction' \cite{Ginzburg}, which defines a non-positively graded $dg$-algebra $\scrJ(Q,W)$ with $H^0(\scrJ(Q,W)) = J(Q,W)$ the Jacobi algebra (i.e., the path algebra modulo cyclic derivatives of the potential).
We always assume that $W$ has no quadratic terms; then the category $\scrC$ has a spherical object for each vertex of $Q$, and $W$ encodes the non-trivial products in a cyclic $A_{\infty}$-structure on the total endomorphism algebra of those objects (this is the Koszul dual of the Ginzburg algebra). By a `symplectic geometric model' for $\scrC(Q,W)$ we mean a configuration of Lagrangian spheres $\{L_v\}_{v \in \mathrm{Vert}(Q)} \subset (X,\omega)$ such that the number of edges between $v$ and $v'$ gives the rank of $HF^*(L_v,L_{v'})$ and the endomorphism algebra $\bigoplus_{v,v'} HF^*(L_v,L_{v'})$ is $A_{\infty}$-quasi-isomorphic to that determined by the potential $W$. Such models are especially useful when the $\{L_v\}$ (split-)generate the compact category $\scrF(X,\omega)$. One can also seek algebro-geometric models for $\scrC(Q,W)$, i.e., embeddings $\scrC(Q,W) \hookrightarrow \scrD(Y)$ into the derived category of some quasi-projective $Y$ (ideally with image some naturally characterised subcategory), or indeed algebraic models in terms of realisations in derived categories of modules over a non-commutative algebra. All of these will play a role in the sequel.

\section[Type A Milnor fibres]{Type $\boldsymbol{A}$ Milnor fibres}
\label{sec:type-A-Milnor-fibres}

We consider the $A^{(n)}_k$-Milnor fibre of complex dimension $n$, that is the affine hypersurface
\begin{equation} \label{eqn:type A}
 \bigl\{ x_1^2+\dots +x_n^2 + p_{k+1}(z) = 0 \bigr\} \subset \bigl(\bC^{n+1}, \omega_{\rm st}\bigr)
\end{equation}
with $p_{k+1}$ a degree $k+1$ polynomial with distinct zeroes. This is symplectically a~plumbing of~$k$ copies of $T^*S^n$, associated to the $A_k$-quiver:

\begin{figure}[h]
	\centering
	\begin{tikzpicture}
	\node[circle,draw, fill, minimum size = 2pt,inner sep=1pt] at (0,0) {};
	\node[circle,draw, fill, minimum size = 2pt,inner sep=1pt] at (8,0) {};
	\node[circle,draw, fill, minimum size = 2pt,inner sep=1pt] at (1.5,0) {};
	\node[circle,draw, fill, minimum size = 2pt,inner sep=1pt] at (3,0) {};
	
	\draw[->,shorten >=8pt, shorten <=8pt] (0,0) to (1.5,0);
	\draw[->,shorten >=8pt, shorten <=8pt] (1.5,0) to (3,0);
	
	\path (3,0) to node {\dots} (8,0);
	\node [shape=circle,minimum size=2pt, inner sep=1pt] at (4.5,0) {};
	\draw[->,shorten >=8pt, shorten <=8pt] (3,0) to (4.5,0);
	
	\node [shape=circle,minimum size=2pt, inner sep=1pt] at (6.5,0) {};
	\draw[->,shorten >=8pt, shorten <=8pt] (6.5,0) to (8,0);
	\end{tikzpicture}
	\caption{The $A_k$-quiver.}
\end{figure}

It comes in a family over the configuration space $\mathrm{Conf}_{k+1}(\bC)$ parametrizing polynomials of degree $k+1$ with distinct zeroes, which yields a representation
\[
\rho\colon \ {\rm Br}_{k+1} \to \pi_0\Symp_{\rm ct}\bigl(A_k^{(n)}\bigr)
\]
from the classical braid group on $k+1$ strings. Combining the work of many authors one has the following assertion.

\begin{prop}[{\cite[Corollary 1.4]{Khovanov-Seidel}}]
The representation $\rho$ is injective for every $k$, $n$.
\end{prop}

\begin{prop}[{\cite[Theorem~1.4]{Evans} and \cite[Theorem~1.2]{Wu}}]
When $n=2$, $\rho$ is an isomorphism.
\end{prop}

\begin{prop}[{\cite[Corollary 17.17]{Seidel:FCPLT} and \cite[Theorem A]{DRE}}]
 For some pairs $k$, $n$ the monodromy $\rho$ is not surjective.
\end{prop}

The first result is based on direct computations of Floer cohomology; the second relies on special features of holomorphic curves in four dimensions; the third result uses `reframings' of twists by exotic diffeomorphisms of spheres (for instance, it applies when $n=8$ and $k\in\{1,2\}$). This last result makes the following interesting statement.

\begin{thm}\label{thm:milnor case}
For $n\geq 3$ and arbitrary $k$, the map $\rho$ splits.
\end{thm}

We will assume $n\geq 3$ throughout the rest of this section. We first establish the slightly simpler proposition.

\begin{prop}\label{prop:mod center}
There is a map $p\colon \pi_0\Symp_{\rm ct}\bigl(A_k^{(n)}\bigr) \to {\rm Br}_{k+1}/Z({\rm Br}_{k+1})$ for which the composite $p \circ \rho$ is the natural projection.
\end{prop}

\begin{proof}
Since $\pi_0\Symp_{\rm ct}\bigl(A_k^{(n)}\bigr)$ maps canonically to $\pi_0\Symp^{\rm gr}\bigr(A_k^{(n)}\bigr)$, it acts on the space of stability conditions on the compact Fukaya category $\scrF$. This has a distinguished component in the stability space coming from tilts of chambers associated to the standard heart. Qiu and Woolf \cite[Theorem~B]{Qiu-Woolf} prove that in fact the space of stability conditions is connected, so the whole of $\Auteq(\scrF)$ acts on this unique component. The autoequivalences preserving the given component have been computed in \cite[Theorem~9.9]{BS} when $n=3$ and in general in \cite[Theorem~7.16]{Ikeda} and \cite[Theorem~B]{Qiu-Woolf}. The complete answer is as follows. We can identify $\scrF$ with the category of modules over the endomorphism algebra of the $A_k$-chain of Lagrangian `core' spheres. Ikeda proves that
\begin{equation} \label{eqn:full aut group for A_k}
\Auteq\bigl(\scrF\bigl(A_k^{(n)}\bigr)\bigr) = ({\rm Br}_{k+1} \times \bZ) / \bZ,
\end{equation}
where the $\bZ$-factor in the numerator is generated by the shift, and the embedding of the denominator sends the generator to $(\Delta, d_{k,n})$ where $\Delta$ generates the centre and
\begin{equation} \label{eqn:dkn}
d_{k,n} = (k+1)(n-2) + 2.
\end{equation}

There are thus exact sequences
\begin{equation*}
1 \to {\rm Br}_{k+1} \to \Auteq\bigl(\scrF\bigl(A_k^{(n)}\bigr)\bigr) \to \bZ/d_{k,n}\bZ \to 1;
\end{equation*}
and
\begin{equation*}
1 \to \bZ \to \Auteq\bigl(\scrF\bigl(A_k^{(n)}\bigr)\bigr) \to {\rm Br}_{k+1}/Z({\rm Br}_{k+1}) \to 1
\end{equation*}
and the latter admits a map from the sequence
\[
1 \to \bZ \to \pi_0\Symp^{\rm gr}_{\rm ct}\bigl(A_k^{(n)}\bigr) \to \pi_0\Symp_{\rm ct}(A_k^{(n)}) \to 1,
\]
where we are using the fact that the forgetful map $ \pi_0\Symp^{\rm gr}_{\rm ct}\bigl(A_k^{(n)}\bigr) \to \pi_0\Symp_{\rm ct}\bigl(A_k^{(n)}\bigr)$ splits by Lemma~\ref{lem:forgetting-gradings-splits}.
The existence of a map $\pi_0\Symp_{\rm ct}\bigl(A_k^{(n)}\bigr) \to {\rm Br}_{k+1} / Z({\rm Br}_{k+1})$ follows.

The map $\rho$ takes the $i$-th standard generator of the braid group to the Dehn twist in the $i$-th Lagrangian sphere in the $A_k$ chain which is the compact core of $A_k^{(n)}$. Under Ikeda's isomorphism~\eqref{eqn:full aut group for A_k}, the braid group factor is generated by the corresponding spherical twists. That implies that $p \circ \rho$ is the natural projection.
\end{proof}

\begin{rmk}
Since the centre of the braid group acts by a shift on $\scrF(X)$, it acts trivially on any simplicial complex constructed from the set of (ungraded) Lagrangian spheres in $X$. One can detect the centre of the braid group through its action on the ungraded category underlying the partially wrapped Fukaya category studied in \cite[Section~20]{Seidel:FCPLT}, which adjoins a single Lefschetz thimble to the $A_k$-chain of spheres and wraps `infinitesimally'. However, the space of stability conditions on this partially wrapped Fukaya category has not been computed.
\end{rmk}

Via the splitting of the map from graded to ungraded symplectomorphisms over the compactly supported subgroup of Lemma~\ref{lem:forgetting-gradings-splits}, we have a representation
\[
\rho\colon \ \pi_0\Symp_{\rm ct}\bigl(A_k^{(n)}\bigr) \to \Auteq\bigl(\scrF\bigl(A_k^{(n)}\bigr)\bigr)
\]
and to prove Theorem~\ref{thm:milnor case} it remains to show that the composite $\chi \circ \rho\colon \pi_0\Symp_{\rm ct}\bigl(A_k^{(n)}\bigr) \longrightarrow \bZ/d_{k,n}\bZ$ vanishes. This will follow directly from a geometric picture of the target.

 Recall the cluster category $\mathcal{C}_{k,n}$ of type $(k,n)$ is by definition the orbit category of the path category of the underlying quiver of $A_k^{(n)}$, given by dividing out by the composition $\tau^{-1}[n-2]$ of the given power of the shift and the inverse of the Auslander--Reiten translation $\tau$, see, e.g.,~\cite{Baur-Marsh}. It is naturally a ${\rm CY}_{n-1}$-category.

 Recall the integer $d_{k,n}$ introduced in \eqref{eqn:dkn}.
 Let $(S,M)$ be the marked bordered surface which is the disc $S$ with a collection of $d_{k,n}$ boundary marked points, equivalently of a regular $d_{k,n}$-gon.
 We think of the group $\bZ/d_{k,n}\bZ$ as the oriented mapping class group $\Gamma(S, M)$ of $(S,M)$.

 \begin{lem}
 There is a natural homomorphism
 \[
\Auteq(\mathcal{C}_{k,n}) \longrightarrow \Gamma(S,M) = \bZ/d_{k,n}\bZ.
\]
 \end{lem}

 \begin{proof}
 This is a variation on the proof of \cite[Theorem 9.9]{BS}, see also \cite[Theorem~7.15]{Ikeda}, and is exactly the result \cite[Theorem~4.7]{Brustle-YuQiu} in the special case when the cluster category is ${\rm CY}_2$. The paper \cite{Baur-Marsh} gives a combinatorial description of the cluster category. Its indecomposable objects are given by the $n$-diagonals of a $d_{k,n}$-gon (i.e., diagonals which decompose the polygon into an $n$-gon and a remainder) and the cluster tilting objects biject with $n$-angulations of the $d_{k,n}$-gon (maximal sets of disjoint-in-the-interior $n$-diagonals; the diagonals may intersect at the boundary, and any such maximal set has $k$ members). See \cite[Proposition~5.4]{Baur-Marsh} and \cite[Theorem~4.2]{Jacquet-Malo}. (In the special case, where the cluster category is ${\rm CY}_2$, one studies actual triangulations of the polygon.) An illustration with $n=4$, $k=3$ is given in Figure~\ref{fig:n-diagonal}. There is a (directed) `flip' orientation on $n$-angulations: a single $n$-diagonal is a diameter in a unique $(2n-2)$-gon, and one can rotate the ends of the given diagonal one place clockwise in this subpolygon to obtain a new $n$-angulation, cf. the labelled edge $\delta$ contained in the shaded $(2n-2)$-gon in Figure~\ref{fig:n-diagonal} and the forward flip along $\delta$ on the right.

\begin{figure}[ht]\centering

\begin{tikzpicture}[
 scale=2,
 dot/.style={circle,fill=black,inner sep=1.2pt},
 every node/.style={font=\small}
]

\def\r{1}


\begin{scope}

\draw (0,0) circle (\r);

\foreach \i in {1,...,10}{
 \coordinate (L\i) at ({( \i-1)*36}: \r);
}

\fill[gray!20]
 (L3)
 arc[start angle=72,end angle=108,radius=\r]
 -- (L7)
 arc[start angle=216,end angle=324,radius=\r]
 -- cycle;

\draw (L4)--(L7);
\draw (L7)-- node[midway,fill=white,inner sep=1pt] {$\delta$} (L10);
\draw (L10)--(L3);

\foreach \i in {1,...,10}{
 \fill (L\i) circle (0.018);
 \node at ({(\i-1)*36}:1.15) {\i};
}

\end{scope}


\begin{scope}[xshift=3cm]

\draw (0,0) circle (\r);

\foreach \i in {1,...,10}{
 \coordinate (R\i) at ({(\i-1)*36}: \r);
}

\draw (R4)--(R7);
\draw (R4)--(R9);
\draw (R3)--(R10);

\foreach \i in {1,...,10}{
 \fill (R\i) circle (0.018);
 \node at ({(\i-1)*36}:1.15) {\i};
}

\end{scope}
\end{tikzpicture}
\caption{An $n$-angulation $(n=4$, $k=3$, $d_{k,n}=10)$ and its forward flip along the $n$-diagonal $\delta$.}\label{fig:n-diagonal}
\end{figure}

 The cluster exchange graph of $\mathcal{C}_{k,n}$ (with vertices being simple hearts and edges arising from tilts) then has a purely combinatorial description in terms of directed flips of $n$-angulations, cf.~\cite[Theorem~5.6]{Baur-Marsh} or \cite[Theorem~6.3]{Ikeda}.

 The orientation-preserving mapping class group of the marked disc acts (with finite stabilisers) on $n$-angulations and hence on the cluster exchange graph.

 The $A_k$-quiver has an associated ${\rm CY}_n$-category, via the Ginzburg construction, which also has a heart exchange graph. It is a non-trivial theorem of Buan, Reiten and Thomas~\cite{BRT} that one can associate to a heart in the Ginzburg category a `silting collection' in the quotient cluster category, and that this defines a map from the heart exchange graph upstairs to the cluster exchange graph downstairs.
 In \cite[Theorem~8.6]{King-QiuYu-cluster}, the authors show that the braid group acts fibrewise with respect to this map and moreover acts transitively on each fibre. This then establishes all the conditions needed for the argument of \cite[Theorem~9.9]{BS}.
 \end{proof}

\begin{lem}
There is an equivalence $\scrR\bigl(A_k^{(n)}\bigr) \simeq \mathcal{C}_{k,n}$.
\end{lem}

\begin{proof} This is observed \emph{en passant} in the discussion before \cite[Conjecture~1.5]{Lekili-Ueda}. A more detailed identification of the cluster category and Rabinowitz Fukaya category for a Dynkin quiver is given in \cite[Corollary~1.4]{BJK}. \end{proof}

\begin{lem}
 There is a commuting diagram
 \[
\xymatrix{
\Auteq\bigl(\scrF\bigl(A_k^{(n)}\bigr)\bigr) \ar^{\mathrm{Ikeda}}[rr] && \bZ/d_{k,n}\bZ \ar@{=}[d] \\
\Auteq\bigl(\scrW\bigl(A_k^{(n)}\bigr)\bigr) \ar[u] \ar[r] & \Auteq\bigl(\scrR\bigl(A_k^{(n)}\bigr)\bigr) \ar[r] & \bZ/d_{k,n}\bZ.
}
 \]
\end{lem}

\begin{proof}
On any Weinstein manifold, the Rabinowitz Fukaya category is algebraically determined from the wrapped category as the `categorical formal punctured neighbourhood of infinity' (given by quotienting out proper modules), cf.~\cite{GGV} and the definition in Section~\ref{Sec:flavours}, which defines the first arrow on the bottom line. The wrapped and compact Fukaya categories are in this case Koszul dual,\footnote{Since $n\geq 3$, the compact category $\scrF(A_k^{(n)})$ is a Koszul dualising subcategory of $\scrW(A_k^{(n)})$; this follows by combining results of~\cite{AbSm-plumbing} and~\cite{Kalck-Yang}, cf.~the discussion of Koszul duality later in Section~\ref{sec:marked-bordered-surfaces} (which encompasses this case when $n=3$, viewing the $A_k^{(n)}$-Milnor fibres as special cases of the 3-folds associated to marked bordered surfaces; the arguments when $n>3$ are identical).} which defines the left vertical arrow.

Commutativity of the diagram amounts to the fact that both horizontal arrows are given by the action of autoequivalences on the heart exchange graph. In more detail, Ikeda \cite[Section~7]{Ikeda} constructs the top arrow directly using the action on the heart exchange graph. The quotient of the heart exchange graph by the braid group action is identified with the cluster exchange graph of $\scrR\bigl(A_k^{(n)}\bigr)$ by \cite{King-QiuYu-cluster}. Furthermore, the cluster exchange graph is purely combinatorial, as discussed above, and the bottom right map is constructed using the action on the cluster exchange graph.
\end{proof}

\begin{rmk}
 It seems likely that the simplicial automorphism group of the cluster exchange graph is given by the mapping class group of the marked surface, which in this case is just the given cyclic group. This is known in the case of triangulations (here corresponding to the 3-fold case), but does not seem to have been written down for $m$-angulations in general. That would give a more direct argument for the previous result.
\end{rmk}

Since compactly supported symplectomorphisms of $X$ act on both $\scrF$ and $\scrW$ and act trivially on the quotient category $\scrR(X)$, we immediately infer that the composite map
\[
\pi_0\Symp_{\rm ct}(X) \to \Auteq\bigl(\scrF\bigl(A_k^{(n)}\bigr)\bigr) \to \bZ/d_{k,n}\bZ
\]
vanishes, hence the first arrow lifts to ${\rm Br}_{k+1}$, as required.

The case of type $A$ Milnor fibres raises an obvious question.

\begin{qn}\label{qn:general-Milnor}
 Let $X$ be the Milnor fibre of an isolated hypersurface singularity of dimension $\geq 2$. Let $\mathcal{M}_X$ denote the universal unfolding of the singularity. Does symplectic monodromy split over $\mathcal{M}_X$?
\end{qn}

\begin{rmk}
 The dimension hypothesis is necessary. Indeed, already in type $A$ in the lowest-dimensional case $n=1$ in \eqref{eqn:type A}, the monodromy representation defines an embedding of the braid group into the mapping class group of a punctured surface. Such embeddings do not split (for instance, the mapping class group typically has trivial abelianisation, and the braid group has abelianisation $\bZ$; there are many more subtle obstructions, e.g.,~the induced map on stable homology is trivial in positive degrees, cf.~\cite{Song-Tillmann}).
\end{rmk}

\section{Quivers and threefolds from marked bordered surfaces}\label{sec:marked-bordered-surfaces}

The type $A$ Milnor fibres $A_k^{(n)}$ have a quite different generalisation from that of Question \ref{qn:general-Milnor}, namely affine Calabi--Yau varieties which are Lefschetz fibrations with general fibre $T^*S^{n-1}$ over a non-compact surface.
We will focus on the case of $3$-folds fibred in $T^*S^2$. In this case, the space of stability conditions is not known to be connected, but we can access an analogue of Proposition~\ref{prop:mod center} and Theorem~\ref{thm:milnor case} using knowledge of spherical objects. (Note that in more complicated cases the relevant braid-type group may have no center, in which case these two results are essentially equivalent.)

Let $S$ be an oriented surface of genus $g$ with a non-empty set of $d$ marked points $\{p_1,\dots,p_d\}$. Fix an integer $d_i > 2$ for each $p_i$; we encode this information by a vector $\bf{d}$. We take the real blow-up of $S$ at the $\{p_i\}$ and mark $d_i-2$ points on the $i$-th boundary component; let $M$ denote the resulting collection of marked points. We consider ideal triangulations of $S$ in the sense of \cite{FST}, i.e., maximal collections of disjoint arcs, where an arc is a smooth embedded path connecting points of $M \subset \partial S$, with interior in $S \backslash \partial S$. (To get something visibly resembling a triangulation one should include arcs in $\partial S$ joining points of $M$, but these are \emph{not} included in the set of arcs.) A triangle is `interior' if its closure meets $\partial S$ only in $M$. Any ideal triangulation has
\begin{equation} \label{eqn:edges and faces}
6g-6+ \sum_i (d_i+1) \ \mathrm{edges} \qquad \mathrm{and} \qquad 4g-4+\sum_i d_i \ \mathrm{faces}
\end{equation}
(we are in the case with `no punctures' in the terminology of \cite{BS}, and we are not including the boundary arcs between points of $M$ as edges).

We will be interested in the cellulation which is dual to the given ideal triangulation, whose vertices will correspond to simple zeroes of a quadratic differential below.

\begin{figure}[ht]
\centering
 \includegraphics[scale=0.4]{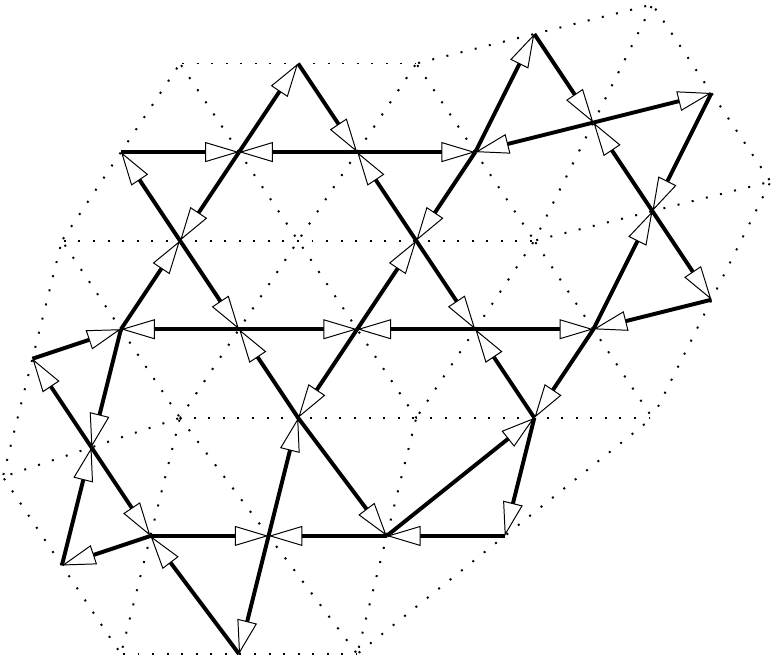}
 \caption{Quiver inscribed in a triangulation.}\label{fig:quiver}
\end{figure}

There is a quiver associated to the ideal triangulation, with vertices the mid-points of the edges and an arrow when two open arcs share a vertex (in the direction in which one is anticlockwise from the other). See Figure~\ref{fig:quiver} for an illustration\footnote{This figure is borrowed from \cite{BS}.} (in a region interior to the surface, so $\partial S$ is not visible in the figure).
The resulting collection of arrows in particular inscribes a~clockwise 3-cycle in each interior triangle, so comes with an obvious collection of $3$-cycles: the triangles just inscribed. This defines a canonical potential
\[
W =\sum (\Circlearrowright \, \mathrm{triangles})
\]
and hence via the Ginzburg construction (cf.~the brief recollection at the end of Section~\ref{Sec:flavours}) a ${\rm CY}_3$-category $\scrC(Q,W) $. One can show this category depends only on the pair $(S,M)$ or equivalently on $(g,\bf{d})$. More precisely, flips of the triangulation induce mutations of quivers-with-potential, which yield quasi-equivalent categories by \cite{Keller-Yang}. We may write $\scrC(S,M)$ or $\scrC(g,\mathbf{d})$ for a category in this equivalence class.

Let $\phi$ be a meromorphic quadratic differential on a surface $S$ with poles of order $d_i$ at $p_i$ and with distinct zeroes. Let $D\subset S$ be the set of poles. Let $Y = Y_{\phi} \to S \setminus D$ be the Lefschetz fibration with smooth fibre $T^*S^2$ and with nodal singular fibres over the zeroes of $\phi$. This is a non-compact affine Calabi--Yau 3-fold, which moves in a family over the complex orbifold $\mathrm{Quad}(g,\mathbf{d})$ of quadratic differentials with the prescribed pole orders and with simple zeroes. An embedded `matching path' in $S$ between zeroes of $\phi$ defines a Lagrangian $3$-sphere in $Y_{\phi}$.

\begin{lem}
An ideal triangulation defines a collection of matching paths, which form the edges of the dual cellulation.
\end{lem}

See \cite{Smith:quiver} for a proof, and further background.

 Let $\Gamma(S,M)$ denote the mapping class group of the associated marked bordered surface.\footnote{For comparison to \cite{BS,Smith:quiver}, note that we are in the case of an unpunctured surface with non-empty boundary.} Let $\mathrm{Sph}(S,M)$ denote the group generated by spherical twists in matching paths dual to the edges of the ideal triangulation.

\begin{thm}[{\cite{BS, Smith:quiver}}]\label{Thm:BS}
There is a distinguished connected component $\Stab^{\dagger}(\scrC(S,M))$ of the space of stability conditions on $\scrC(S,M)$, and a subgroup $\Auteq^{\dagger}(\scrC) \subset \Auteq(\scrC)$ of autoequivalences preserving that component, with the following properties:
\begin{itemize}\itemsep=0pt
\item $\Stab^\dagger(\scrC(S,M))/\Auteq^\dagger(\scrC) \cong \Quad(g,\mathbf{d})$ is the complex orbifold whose points parametrise biholomorphism classes of pairs of a Riemann surface $S$ of genus $g$ and a meromorphic quadratic differential on $S$ with poles of orders $d_i>2$ and with simple zeroes.
\item There is an exact sequence
\begin{equation*} 
1 \to \mathrm{Sph}(S,M) \to \Auteq^{\dagger}(\scrC(S,M)) \to \Gamma(S,M) \to 1.
\end{equation*}
\item There is an embedding $\scrC(Q,W) \hookrightarrow \scrF(Y_{\phi})$ for $(Q,W)$ the quiver with potential associated to any ideal triangulation of $(S,M)$.
\end{itemize}
\end{thm}

\begin{rmk}
Let $Z \subset \mathrm{interior}(S) \subset S\setminus M$ denote a set of points in bijection with the (assumed simple) zeroes of a quadratic differential $\phi$ of the prescribed polar type. We can consider the mapping class group $\Gamma(S,M;Z)$ of mapping classes which fix $M$ setwise and which preserve the finite set $Z$ setwise. This fits into a short exact sequence
\[
1 \to {\rm SBr}(S,Z) \to \Gamma(S,M;Z) \to \Gamma(S,M) \to 1,
\]
where the first term is the `surface braid group' of $(S,Z)$, i.e.,~the fundamental group of the space of configurations of unordered $|Z|$-tuples of points on~$S$. The surface braid group has two kinds of generators: braid-type maps associated to a matching path between two distinct points of~$Z$, and point-pushing maps associated to a non-trivial loop on $S$ and a single point of~$Z$. The paper~\cite{QiuYu:spherical-twist-group} has shown that $\mathrm{Sph}(S,M)$ is the subgroup of ${\rm SBr}(S,Z)$ generated by the first type of map.

Let $\Sigma \to S$ be the spectral curve, which is the double cover branched at $Z \subset S$, and let $\iota\colon \Sigma \to \Sigma$ denote the covering involution. Let $\Gamma(\Sigma,\iota)$ denote the hyperelliptic-type mapping class group of isotopy classes of diffeomorphisms of $\Sigma$ which commute with $\iota$. Associating to a~braid-twist in a matching path $\gamma$ on $S$ with end-points in $Z$ the Dehn twist in the $\iota$-invariant simple closed curve $\hat\gamma \subset \Sigma$ which is the lift of $\gamma$ shows that there is an embedding
\[
\mathrm{Sph}(S,M) \hookrightarrow \Gamma(\Sigma,\iota).
\]
According to the work of \cite{Birman-Hilden} on symmetric mapping class groups, cf.~also \cite{Margalit-Winarski}, and noting that we might as well restrict to the part of the autoequivalence group which maps to $\Gamma(S,M_{=})$ which fixes $M$ and hence $\partial S$ pointwise, we find that when $\chi(S)<0$, then $\mathrm{Sph}(S,M)$ is characterised as the centraliser of the involution $\iota$ in $\Gamma(\Sigma,\partial\Sigma)$. (We do not need to divide by $\iota$ since the mapping class group of a surface with non-empty boundary, fixing the boundary pointwise, is torsion-free, so $\iota$ does not define an element of $\Gamma(\Sigma,\partial\Sigma)$.)
\end{rmk}

The group $\mathrm{Sph}(S,M)$ acts via spherical twists which are geometrically realised by Dehn twists in Lagrangian matching spheres, so there is a natural representation
\[
\mathrm{Sph}(S,M) \to \pi_0\Symp_{\rm ct}(Y). \]

The following generalises (the 3-dimensional instance of) Proposition~\ref{prop:mod center}.

\begin{thm} \label{thm:bordered surface case}
The monodromy map $\rho\colon \mathrm{Sph}(S,M) \to \pi_0\Symp_{\rm ct}(Y)$ splits.
\end{thm}

Recall that a spherical object is reachable if it is one of the simples in a reachable heart, i.e.,~one obtained from the standard heart by a finite sequence of tilts.

\begin{lem}
The set of reachable spherical objects in $\scrC(Q,W)$ is in bijection with the set of Lagrangian matching spheres.
\end{lem}

\begin{proof} This is proved in~\cite{King-QiuYu} (note we are in the `unpunctured case' which they treat). \end{proof}

\begin{lem} \label{lem:quis to matching}
Every Lagrangian $3$-sphere in $Y_{\phi}$ is quasi-isomorphic to the Lagrangian $3$-sphere associated to some $($unique up to isotopy$)$ embedded matching path.
\end{lem}

\begin{proof}
Let $L \subset Y_{\phi}$ be a Lagrangian $3$-sphere. Then $L$ lifts to a Lagrangian $\tilde{L}$ in the universal cover $\pi\colon \tilde{Y} \to Y$, which is a Lefschetz fibration $\tilde{Y} \to \bC$ with fibre $T^*S^2$ and infinitely many Lefschetz singular fibres.
This is an `infinite' Milnor fibre of type $A$ in dimension~3.
Being compact, any lift of $L$ lives in a finite subspace which is the total space of a~fibration \smash{$A_k^{(3)}$} for some $k\gg 0$.
Spherical objects in \smash{$A_k^{(3)}$} are classified; for instance, this follows from the result of~\cite{King-QiuYu} together with the fact that the heart exchange graph is connected in type~$A$ \cite{Ikeda}, which means that all hearts are reachable.
It follows that there is a unique embedded matching path~$\gamma$, with associated Lagrangian sphere $\tilde{L}_\gamma$, such that~$\tilde{L}$ is quasi-isomorphic to $\tilde{L}_\gamma$ in the Fukaya category of~\smash{$A_k^{(3)}$}.

Let $\sigma$ denote the immersed matching path downstairs which is the image of $\gamma$, and let $L_\sigma$ be the corresponding immersed matching sphere, which is the image of $\tilde{L}_\gamma$ under $\pi$. We want to show that $L$ and $L_\sigma$ are quasi-isomorphic and that $\sigma$ is embedded. In order to do this, we will construct a quasi-isomorphism upstairs and descend it.

Let $\tilde{L}_\gamma'$ be a small Hamiltonian perturbation of $\tilde{L}_\gamma$.
By a small perturbation of $L$ pulled back to $\tilde{L}$, we can assume that the pairs $(\tilde{L}, \tilde{L}_\gamma )$
and $(\tilde{L}, \tilde{L}_\gamma')$ both intersect transversally.
Let $L_\sigma'$ be the image under $\pi$ of $\tilde{L}_\gamma'$.
Note that $\gamma$ being embedded means that $\sigma$ bounds no teardrops disjoint from the marked points, which implies the associated immersed (cleanly self-intersecting) Lagrangian 3-sphere $L_{\sigma} \subset Y_{\phi}$ is unobstructed. Indeed, the obstruction term $\mathfrak{m}_0$ in immersed Floer theory arises from holomorphic discs with a unique boundary puncture; exactness of the immersion means that this puncture must map to one of the transverse self-intersections, defining a teardrop disc, cf.~\cite[Section 13.3]{Akaho-Joyce} for a discussion of related examples.

Pick a regular $J$ on $ Y_{\phi}$ and pull it back to $\tilde{J}$ on $\tilde{Y}$. Since $\tilde{Y} \to Y$ is an unramified cover, and the domain of a holomorphic disc $u\colon D \to Y$ is simply-connected, it admits a unique lift to $\tilde{Y}$ given any fixed point in the preimage, and since the almost complex structure on $\tilde{Y}$ is pulled back from~$Y$, this lift is again holomorphic. We will apply this where $D$ is a holomorphic strip or triangle in~$Y$, and the point lift which fixes a preimage is an intersection point (using that the map extends over the boundary punctures). This lifting property of holomorphic strips therefore means that the projection map $\pi$ defines a chain map
\begin{equation} \label{eq:chain-map-for-covering-for-Lag-S^3s}
 {CF}_{\tilde{Y}} \bigl(\tilde{L}_\gamma, \tilde{L}; \tilde{J}\bigr) \to {CF}_Y (L_\sigma, L; J)
\end{equation}
and similarly for $\bigl(\tilde{L}, \tilde{L}_\gamma\bigr) $ and $\bigl(\tilde{L}_\gamma, \tilde{L}_\gamma'\bigr)$.

Pick cycles $\alpha \in {CF}_{\tilde{Y}} \bigl(\tilde{L}, \tilde{L}_\gamma' \bigr) $, $\beta \in {CF}_{\tilde{Y}} \bigl(\tilde{L}_\gamma, \tilde{L}\bigr) $ such that
\begin{equation*}
 \mu_2 (\alpha, \beta) = e \in {CF}_{\tilde{Y}} \bigl(\tilde{L}_\gamma, \tilde{L}_\gamma'\bigr),
\end{equation*}
where $e$ represents the cohomological unit. The images of $\alpha$, $\beta$ and $e$ under $\pi$ all represent cycles by~\eqref{eq:chain-map-for-covering-for-Lag-S^3s}, say $\pi(\alpha)$, $\pi(\beta)$ and $\pi(e)$. By the same unique lifting for based holomorphic triangles,
there is a non-trivial Floer product
\begin{equation*}
 \mu_2 (\pi(\alpha), \pi(\beta)) = \pi(e) \in {CF}_{Y} (L_\sigma, L_\sigma')
\end{equation*}
and $\pi(e)$ again represents the cohomological unit.

Taking the product in the other order (and using $L$, $L_\sigma$ and a small Hamiltonian perturbation of~$L$), we deduce that $L$ and $L_\sigma$ are quasi-isomorphic.

Finally, we can assume that $\sigma$ self-intersects minimally. There is a fibrewise involution on the Lefschetz fibration \smash{$A_k^{(3)}$}, with fixed locus a two-dimensional Riemann surface of finite type, which descends to a fibrewise involution on a compact subdomain of~$Y_{\phi}$ containing $L_{\sigma}$. The matching sphere $\tilde{L}_{\gamma}$, and hence its image $L_{\sigma}$, can be constructed to be setwise preserved by this involution. If $\sigma$ is not embedded, then by the arguments of \cite[Section~5]{Khovanov-Seidel} (working with a $\bZ/2$-equivariant $J$ and achieving equivariant transversality since there are no curves in the fixed locus for topological reasons), $HF(L_{\sigma},L_{\sigma};\bZ/2)$ has rank $>2$, which is a contradiction since $L_{\sigma}$ is spherical.
\end{proof}

\begin{cor}
 $\scrC(Q,W)$ is a symplectomorphism-invariant subcategory of $\scrF(Y_{\phi})$, and all Lagrangian spheres are reachable spherical objects.
\end{cor}

We now know that $\pi_0\Symp_{\rm ct}$ preserves reachable spherical objects, but not at this point reachable hearts.

\begin{lem} \label{lem:HF rank one}
Let $L$ and $L'$ be two Lagrangian spheres in $Y_\phi$. Then $HF^*(L,L')$ has rank one exactly when the two corresponding matching paths share one end-point.
\end{lem}

\begin{proof}
Let $\gamma$ and $\gamma'$ denote the relevant matching paths. By isotopy of these (which induces a~Hamiltonian isotopy of the corresponding matching Lagrangian 3-spheres), we can assume that the two paths meet minimally. If two sub-arcs of~$\gamma \cup \gamma'$ bound a loop which is non-trivial in~$\pi_1(S)$ then there can clearly be no Floer strips with those sub-arcs as boundary. Suppose two sub-arcs of the paths bound a region which contains a Lefschetz critical point of the fibration. Then there are no holomorphic strips which project to that region; this follows from the same vanishing result that underlies the proof of the long exact triangle in Floer theory (cf.~\cite[Lemma 17.4, Remark~17.8]{Seidel:FCPLT} and the argument `bubbling off' a standard Lefschetz fibration in \cite[Section 17d]{Seidel:FCPLT}). Given this, the result is straightforward, and closely modelled on \cite{Khovanov-Seidel} as in the final part of the proof of Lemma~\ref{lem:quis to matching}.
\end{proof}

\begin{lem} \label{lem:HF rank two}
Let $L$ and $L'$ be two Lagrangian spheres in $Y_\phi$ which are not quasi-isomorphic. Then $HF^*(L,L')$ has rank two exactly when
\begin{enumerate}\itemsep=0pt
\item[$(1)$] either the two associated matching paths share both endpoints and are otherwise disjoint;
\item[$(2)$] or the two associated matching paths meet transversely at exactly one interior point.
\end{enumerate}
\end{lem}

\begin{proof}
Since $L \not\sim L'$ the associated matching paths $\gamma$ and $\gamma'$ represent different relative homotopy classes in $\pi_1(S,\mathrm{Zeroes})$. If they meet only at both end-points, then since the paths are not homotopic their union traces a non-trivial loop, which means that the union $L_{\gamma} \cup L_{\gamma'}$ bounds no rigid holomorphic strips, so $HF(L,L')$ is indeed rank two. If they meet transversely at an interior point, then $L_{\gamma}$ and $L_{\gamma'}$ meet cleanly along an $S^2$ lying above that point. Putting the matching paths in minimal position and arguing as in the proof of Lemma~\ref{lem:HF rank one} shows that any such $S^2$-component of the intersection locus contributes a two-dimensional vector space to Floer cohomology.
\end{proof}

\begin{lem}\label{lem:two cases}
Suppose $L$, $L'$ are not quasi-isomorphic and $HF^*(L,L')$ has rank two. The two cases of Lemma~{\rm \ref{lem:HF rank two}} cannot be related by a global symplectomorphism.
\end{lem}

\begin{proof}
Assume first that our differentials $\phi$ have at least five zeroes. We claim that the two cases are distinguished by whether or not there is a third Lagrangian sphere $L''$ with $HF^*(L,L'')$ of rank one and $HF^*(L',L'')$ of rank zero. To see this, suppose
$p \in \gamma$ is an end-point of $\gamma$ and $p\not\in \gamma'$. Then we are necessarily in the second case, and the union $\gamma \cup \gamma'$ lies in a disc in $S$ and in particular does not separate. Then one can take a third matching path $\sigma$ which meets $\gamma$ at $p$ but is disjoint from $\gamma'$. Suppose then $\phi$ has exactly four zeroes (so the original ideal triangulation had four faces). By \eqref{eqn:edges and faces}, $S$ is a disc, annulus or a torus with one boundary component. The first two cases correspond to the $A_4$ or affine $A_4$ spaces, where the result is well-known (cf.~the proof of Lemma~\ref{lem:HF rank one}). In the final case, either the two matching paths sweep a non-trivial loop in $\pi_1(T^2)$, and the result is straightforward, or the whole configuration lies in a disc or annulus, reducing one to a previously considered case.

This shows that a suitable $L''$ exists in case (2) of Lemma~\ref{lem:HF rank two}. We claim it cannot exist in case (1). Indeed, if $HF^*(L,L'')$ has rank one, then by Lemma~\ref{lem:HF rank one} $L''$ is associated to a~matching path which shares exactly one end-point in common with $L$ and is otherwise disjoint from it. But then~$L''$ also shares exactly one end-point with the matching path defining $L'$. This means that (after perturbing to achieve transversality) the rank of ${CF}(L',L'')$ is odd, since any interior intersection points of the corresponding matching paths contribute two generators from Morsifying an $S^2$-clean intersection. But then $HF(L',L'') \neq 0$.
\end{proof}

The set $\mathrm{Tri}(S,M)$ of ideal triangulations of $(S,M)$ forms a graph, where two triangulations are connected by an edge if they are related by a flip. (Since we have no punctures, there are no self-folded triangles.)

\begin{lem}
$\pi_0\Symp_{\rm ct}(Y_{\phi})$ acts on the graph of ideal triangulations of $(S,M)$.
\end{lem}

\begin{proof}
Let $f \in \pi_0\Symp_{\rm ct}(Y_{\phi})$. Fix an ideal triangulation $\Delta$, which defines a collection of~matching spheres $\{L_e\}$ indexed by the edges of $\Delta$. We know that for each $e$ the Lagrangian sphere $f(L_e)$ is associated to some new matching path $\gamma_f(e)$. Moreover, Lemmas~\ref{lem:HF rank one}, \ref{lem:HF rank two} and~\ref{lem:two cases} imply that $\gamma_f(e)$ and $\gamma_f(e')$ share exactly one respectively two end-points exactly when $e$ and $e'$ share exactly one respectively two end-points. The collection $\{\gamma_f(e)\}_{e\in \Delta}$ form a~configuration of $4g-4+\sum k_i$ matching paths with no interior intersections, so these must form a maximal set of such arcs, hence are dual to the edges of some ideal trian\-gula\-tion.
\end{proof}

\begin{cor}
$\pi_0\Symp_{\rm ct}(Y_{\phi})$ preserves $\mathrm{Stab}^{\dagger}(\scrC(S,M))$.
\end{cor}

\begin{proof}
It acts on the graph $\mathrm{Tri}(S,M)$, which is connected. This graph is identified \cite{King-QiuYu} with the cluster exchange graph of reachable hearts. Therefore the group sends reachable hearts to reachable hearts, and hence preserves the distinguished component of the space of stability conditions.
\end{proof}

If we fix a generic ideal triangulation with no self-folded triangles, we get a finite collection of Lagrangian 3-spheres $L_e$ indexed by the dual 1-cycles to the edges $e$ of the triangulation, and these meet pairwise transversely. Let $Y_{\rm pl}$ denote the plumbing of copies of $T^*S^3$ according to the graph defined by the dual cellulation to the ideal triangulation. A Stein subdomain of the compact core of $Y_{\rm pl}$ embeds into $Y_{\phi}$, inducing an equivalence on compact and wrapped Fukaya categories. Let $T_e^*$ denote the dual cocore to $L_e$ (so $T_e^*$ is disjoint from $L_{e'}$ for $e' \neq e$, i.e.,~we are not using Lefschetz thimbles but ascending manifolds of a suitable plurisubharmonic function, or cotangent fibres in a `plumbing model').

\begin{lem}\samepage
In the previous situation:
\begin{enumerate}\itemsep=0pt
\item[$(1)$] $Y_{\phi}$ admits a grading structure for which the algebra $\bigoplus_{e,e'} HF^*(L_e,L_{e'})$ is concentrated in non-negative degree and $\bigoplus_{e,e'} HW^*(T^*_e, T^*_{e'})$ is concentrated in non-positive degree.
\item[$(2)$] The inclusion $\scrC(Q,W) \hookrightarrow \scrF(Y_{\phi})$ is an equivalence, i.e.,~the compact cores generate the compact category. Thus $\scrF(Y_{\phi})$ is equivalent to the $dg$-category of finite-dimensional modules over the Ginzburg algebra $\scrJ(Q,W)$.
\item[$(3)$] The compact category $\scrF(Y_\phi)$ and wrapped category $\scrW(Y_{\phi})$ are Koszul dual, and $\scrW(Y_{\phi})$ is equivalent to the $dg$-category of all modules over $\scrJ(Q,W)$.
\end{enumerate}
\end{lem}

\begin{proof}
The first two results follow from \cite{AbSm-plumbing}. Finite-dimensionality of
\[
H^0(\scrW) = \bigoplus_{e,e'} HW^0(T^*_e,T^*_{e'})
\]
follows from the explicit computation in \cite{Karabas-Lee} using that around any puncture, the sum of the corresponding $d_e$ is strictly negative, because the underlying quadratic differential has a pole of order $d_i > 2$. Given that, the argument in \cite[Proposition~4.2 and Lemma~4.6]{AbSm-plumbing} shows that the compact cores generate the compact category. This argument crucially uses that, since the Floer algebra of the cotangent fibres is concentrated in non-positive degree, the ascending degree filtration on any module over $\bigoplus_{e,e'} HF^*(T^*_e, T^*_{e'})$ is a filtration by $A_{\infty}$-submodules. We note \cite[Lemma~4.12]{AbSm-plumbing}, which is a step towards \cite[Proposition~4.2]{AbSm-plumbing}, uses a grading argument which is entirely local near each intersection point, providing one chooses gradings for all core components which are as symmetric as possible, cf. \cite[equation~(4.19)]{AbSm-plumbing}.

Koszul duality then follows from work of Kalck and Yang: finite-dimensionality of $H^0(\scrW)$ (together with concentration of the whole algebra in non-positive degree) and \cite[Proposition~2.5]{Kalck-Yang} shows that $HW^i(\scrW)$ is finite-dimensional for each $i$, which fulfils the hypothesis of \cite[Theorems~2.7 and~2.8]{Kalck-Yang}. The equivalence to the given Ginzburg algebras is established in \cite{Smith:quiver} for the compact category and \cite{Christ, Karabas-Lee} for the wrapped category.
\end{proof}

\begin{lem}
The Rabinowitz Fukaya category of $Y_{\phi}$ is the cluster category $\mathrm{Clust}(Q,W)$ of the quiver with potential $(Q,W)$. Moreover, the cluster tilting objects are in bijection with ideal triangulations.
\end{lem}

\begin{proof}
The first statement follows from the previous result, and the definition of the cluster category for a quiver with potential with finite-dimensional Jacobi algebra, cf.~\cite{Amiot}. In the case of a marked bordered surface, the cluster-tilting objects are classified in \cite{Brustle-Zhang}.
\end{proof}

\begin{lem}
There is a natural surjection $\Auteq(\mathrm{Clust}(Q,W)) \longrightarrow \Gamma(S,M)$
\end{lem}

\begin{proof}
We have an identification between the triangulation graph $\mathrm{Tri}(S,M)$ and the cluster exchange graph, so the given autoequivalence group acts by simplicial automorphisms of the triangulation graph. The group of such simplicial automorphisms was computed in \cite{Disarlo} and is exactly $\Gamma(S,M)$.
\end{proof}

We now have a map
\[
\mathrm{Sph}(S,M) \to \pi_0\Symp_{\rm ct}(Y_{\phi}) \to \Auteq^{\dagger}(\scrC(Q,W))).
\]
Our main result is completed on showing the following.

\begin{lem}
The composite map
\[
\pi_0\Symp_{\rm ct}(Y_{\phi}) \to \Auteq^{\dagger}(\scrC(Q,W)) \to \Gamma(S,M)
\]
vanishes.
\end{lem}

\begin{proof}
The map $\Auteq^{\dagger}(\scrC(Q,W)) \to \Gamma(S,M)$ factors through autoequivalences of the cluster category, hence the composite map factors through the autoequivalence group of the Rabinowitz category, and so the result follows from Lemma~\ref{lem:trivial on Rabinowitz}.
\end{proof}

\section{The double bubble and multi-bubbles} \label{sec:double-bubble}

Recall from Section~\ref{Sec:flavours} that one can associate a 3-dimensional\footnote{There is a generalisation to a `graded quiver with potential' which allows one to associate a ${\rm CY}_d$-category for other values of $d$, cf. \cite{JieRen}.} Calabi--Yau category to a quiver with potential, i.e., a quiver equipped with a cyclic word on the edges.
Arguably the simplest quiver admitting a non-trivial potential is the oriented 2-cycle:
\[
\xymatrix@C=1em{\bullet \ar@/^/[rr]^e && \bullet \ar@/^/[ll]^f }, \qquad \mathrm{potential} \ W = (ef)^2.
\]
(Any potential on this quiver has the form $W = \sum_j a_j (ef)^j = (ef)^k(a_k + \cdots) \simeq (ef)^k$, up to formal change of variable and working over $\bC$.) We write $\scrC$ for the corresponding ${\rm CY}_3$-category.

A symplectic model $X$ for $\scrC$, in the sense of Section~\ref{Sec:flavours}, is given by the affine $(1,1)$-hypersurface in $\bP^2\times\bP^2$. The symplectic completion of this affine variety is also described as the plumbing of two Lagrangian $S^3$s, say $L_1$, $L_2$, along an embedded unknotted circle $S^1$ such that the Morse--Bott Lagrangian surgery $K$ of $L_1$, $L_2$ is again a 3-sphere. Considering suitable pencils and nets of hypersurfaces shows that there is a monodromy representation
\[
{\rm PBr}_3 \to \pi_0\Symp_{\rm ct}(X,\omega)
\]
with the pure braid group generated by the spherical twists $\tau_{L_1}$, $\tau_{L_2}$, $\tau_K$; this can be seen more concretely from the realisation of $X$ as a Morse--Bott--Lefschetz fibration over~$\bC$, which comes in a family of such where one moves the critical values of the fibration. Since the vanishing cycles differ at each singularity, only pure braids show up.

\begin{thm} \label{thm:double bubble case}
The natural monodromy map ${\rm PBr}_3 \to \pi_0\Symp_{\rm ct}(X,\omega)$ splits.
\end{thm}

En route to showing this, it will be helpful to know the following.

\begin{lem} \label{lem:trivial on K-theory}
Any compactly supported symplectomorphism of manifold $X$ acts trivially on $K_{\rm num}(\scrF(X)) \simeq \bZ^2$.
\end{lem}

\begin{proof}
The argument is very close to \cite[Section 3.6]{Keating-Smith}. The natural Euler pairing
\[
K(\scrF(X)) \times K(\scrW(X)) \to \bZ, \qquad (L,\Delta) \mapsto \chi(HF(L,\Delta))
\]
descends to a non-degenerate pairing between $K_{\rm num}(\scrF(X))$ and $K(\scrW(X))$, where $K_{\rm num}$ was introduced in \eqref{eqn:Knum}. The group $K(\scrW(X))$ is generated by the two cotangent fibres; one thus sees that $K_{\rm num}(\scrF(X)) \cong H_3(X;\bZ)$ is generated by the two compact core 3-spheres of the plumbing, and the non-degenerate pairing between $K_{\rm num}(\scrF(X))$ and $K(\scrW(X))$ is identified with the intersection pairing
\begin{equation} \label{eqn:intersection pairing}
H_3(X) \otimes H_3(X,\partial X) \to \bZ.
\end{equation}
The intersection pairing between the two core 3-spheres vanishes, indeed they are smoothly disjoinable (for a more striking and closely related statement see Corollary \ref{cor:disjoinable}). Non-degeneracy of \eqref{eqn:intersection pairing} then implies that $H_3(X) \to H_3(X,\partial X)$ vanishes. But then from the exact sequence
\[
H_3(\partial X) \to H_3(X) \to H_3(X,\partial X)
\]
which is functorial under symplectomorphisms, we see that compactly supported maps (which act trivially on the first term) act trivially on $H_3(X)$.
\end{proof}

There is also a geometric model for $\scrC$ in algebraic geometry. Take
$Y$ to be the crepant resolution $f\colon Y \to \mathrm{Spec}(R)$ of the isolated singularity given by $R = \bC[[u,v,x,y]] / \langle uv-xy(x+y) \rangle$; then $Y$ is a neighbourhood of two floppable $(-1,-1)$-curves meeting once. We consider the subcategory $\scrC_Y$ of the derived category of $Y$ consisting of objects $\scrE$ supported on the exceptional locus and with $Rf_*\scrE = 0$ (this is split-generated by the $\mathcal{O}_{\bP^1}(-1)$ for the two exceptional $\bP^1$'s).

\begin{figure}[ht]\centering
\begin{tikzpicture}[scale=0.5]

\draw[semithick,dashed] (1,0.5) -- (1,2);
\draw[semithick,blue] (1,3) ellipse (0.5cm and 1cm);
\draw[semithick, blue] (1,3.5) arc (150:210:1);
\draw[semithick, blue] (0.9,3.3) arc (35:-30:0.6);

\draw[semithick, red] (1,3) ellipse (0.3cm and 0.8 cm);
\draw[semithick, orange] (1.02,3) arc (200:345:0.25);
\draw[semithick, orange, dotted] (1.02,3) arc (110:70:0.75);

\draw[semithick, red, dashed, rounded corners] [->] (0,2) -- (0,3) -- (0.7,3);
\draw (0,1.5) node {$b$};

\draw[semithick, orange, dashed, rounded corners] [->] (2,2) -- (2,2.5) -- (1.5,3);
\draw (2,1.5) node {$a$};

\draw[fill] (0,0) circle (0.1);
\draw[fill] (2,0) circle (0.1);
\draw[fill] (-2,0) circle (0.1);

\draw[semithick] plot [smooth, tension=1] coordinates { (0,0) (0.5,0.5) (1.25,-0.25) (2,0) };
\draw[semithick] plot [smooth, tension=1] coordinates { (0,0) (-0.5,-0.25) (-1.25,0.25) (-2,0) };

\draw[fill] (7,0) circle (0.1);
\draw[fill] (9,0) circle (0.1);
\draw[fill] (5,0) circle (0.1);

\draw[semithick] plot [smooth, tension=1] coordinates { (7,0) (7.5,0.5) (8.25,-0.25) (9,0) };
\draw[semithick] plot [smooth, tension=1] coordinates { (7,0) (6.5,-0.25) (5.755,0.25) (5,0) };

\draw (6.8,0.4) node {$b$};
\draw (5,0.45) node {$a$};
\draw (9,0.45) node {$a+b$};

\draw[semithick, blue] (6.3,2) -- (7.7,2) -- (7.7,2.2) -- (6.3,2.2) -- (6.3,2);
\draw[semithick, blue] (4.9,1.5) -- (5.1,1.5) -- (5.1,2.5) -- (4.9,2.5) -- (4.9,1.5);
\draw[semithick,blue] (8.5,1.5) -- (8.9,1.5) -- (9.9,2.5) -- (9.5,2.5) -- (8.5,1.5);

\draw[semithick, dashed, gray] (5,1.3) -- (5,0.7);
\draw[semithick, dashed, gray] (7,1.8) -- (7,0.6);
\draw[semithick, dashed, gray] (9,1.3) -- (9,0.7);
\end{tikzpicture}
\caption{Vanishing cycles for Morse--Bott--Lefschetz fibration $X \to \bC$, fibre $(\bC^*)^2$.}\label{Fig:vcycles}
\end{figure}

There are actions of the pure braid group on $X$ and on $D(Y)$ via symplectic parallel transport respectively via flop-flop functors.
Mirror symmetry in this example states the following.

\begin{thm}[{\cite[Proposition 4.27\,(1) and Corollary 5.7]{Smith-Wemyss}}]
There are equivalences
\[
\scrF(X) \simeq \scrC_Y \qquad \mathrm{and} \qquad \scrW(X) \simeq D(Y)/\langle \mathcal{O}_Y \rangle
\]
which entwine the action of the pure braid group.
\end{thm}

\begin{rmk} It seems likely that $D(Y)$ itself is a `partially wrapped' category of $X$; however, it seems hard to characterise $D(Y)$ in a way which lets one pick out (and prove the equivalence to) a distinguished such. Recall that $X$ is the complement of a generic $(1,1)$-divisor in the flag 3-fold, so $X$ admits a projective compactification with smooth boundary divisor. This means there is a Reeb flow on $\partial X$ which rotates the fibres of a circle bundle, and hence a cochain complex computing $SH^*(X)$ which has linear growth at chain level, compare to \cite[Section~4a]{Seidel:bias}. It follows that $SH^0(X)$ has Krull dimension at most~1, and $\operatorname{Spec} SH^0(X)$ has dimension at most~$1$. If $\scrW(X)$ was a categorical crepant resolution of the ring $R$ in the sense of~\cite{Pomerleano}, one would have $\operatorname{Spec} SH^0(X) \cong \operatorname{Spec}R$ being 3-dimensional. Thus, our setting is somewhat different to the one often considered in mirror symmetry for log Calabi--Yau's, as discussed in \cite{Pomerleano}.
\end{rmk}

 The space of stability conditions on $\scrC_Y = \scrC$ was shown to be connected in \cite{Hara-Wemyss}, so the `distinguished' component is preserved by all autoequivalences. They also proved that the `Fourier--Mukai' equivalences (namely those induced from equivalences of $D(Y)$ commuting with $Rf_*$) define a group $\Auteq^{\rm FM}(\scrC) \simeq {\rm PBr}_3$. Unfortunately, the symplectic counterpart of commuting with $Rf_*$ is not immediate. We circumvent this by using the classification of all hearts of bounded $t$-structures on~$\scrC$.

We briefly recall some results on the Deligne groupoid that appears in this setting \cite{Donovan-Wemyss, Hara-Wemyss}: this is a way of labelling a collection of `mutation' equivalences between a finite set of derived categories associated to $R$ by paths in a certain fundamental groupoid, with loops defining autoequivalences.

There is an explicit crepant resolution of $\mathrm{Spec}(R)$ associated to each ordering of the linear factors $\{\ell_1,\ell_2,\ell_3 \} = \{x,y,(x+y)\}$ occurring in its defining equation, given by blowing up the ideals $(u,\ell_i)$ sequentially. The different resolutions arising from different orderings may be isomorphic, but are not canonically isomorphic. Thus, we obtain six different crepant resolutions.

One can place these six resolutions in the chambers of a hyperplane arrangement. General theory \cite[Lemma 3.3]{Smith-Wemyss} identifies these crepant resolutions with `non-commutative crepant resolutions' (NCCRs), which are endomorphism algebras\footnote{An explicit algebraic computation for one of the NCCRs, realising it as a quiver with relations \cite[Lemma~3.5]{Smith-Wemyss}, pins down its $A_{\infty}$-structure, which plays a key role in comparing it to the symplectic geometric models.} of certain maximal Cohen--Macaulay modules over $R$. From this viewpoint, there are six NCCRs associated to the six geometric crepant resolutions, which are related by \emph{mutation functors} $\Phi_1$ and $\Phi_2$, corresponding to mutation between the NCCRs of $R$ at vertex~$1$, and vertex~$2$, respectively.
See \cite[Section~6a]{Smith-Wemyss} for further explanation. These mutation functors form a representation of the Deligne groupoid associated to the given hyperplane arrangement, as follows:
$$
\begin{tikzpicture}[scale=1.3,bend angle=15, looseness=1,>=stealth]
\coordinate (A1) at (135:2cm);
\coordinate (A2) at (-45:2cm);
\coordinate (B1) at (153.435:2cm);
\coordinate (B2) at (-26.565:2cm);
\draw[red!30] (A1) -- (A2);
\draw[black!30] (-2,0)--(2,0);
\draw[black!30] (0,-1.8)--(0,1.8);
\node (C+) at (45:1.5cm) {$\scriptstyle\Db(\Lambda)$};
\node (C1) at (112.5:1.5cm) {$\scriptstyle\Db(\Lambda_1)$};
\node (C2) at (157.5:1.5cm){$\scriptstyle\Db(\Lambda_{12})$};
\node (C-) at (225:1.5cm) {$\scriptstyle\Db(\Lambda_{121})$};
\node (C4) at (-67.5:1.5cm) {$\scriptstyle\Db(\Lambda_{21})$};
\node (C5) at (-22.5:1.5cm) {$\scriptstyle\Db(\Lambda_{2})$};
\draw[->, bend right] (C+) to (C1);
\draw[->, bend right] (C1) to (C+);
\draw[->, bend right] (C1) to (C2);
\draw[->, bend right] (C2) to (C1);
\draw[->, bend right] (C2) to (C-);
\draw[->, bend right] (C-) to (C2);
\draw[<-, bend right] (C+) to (C5);
\draw[<-, bend right] (C5) to (C+);
\draw[<-, bend right] (C5) to (C4);
\draw[<-, bend right] (C4) to (C5);
\draw[<-, bend right] (C4) to (C-);
\draw[<-, bend right] (C-) to (C4);
\node at (78.75:0.9cm) {$\scriptstyle \Phi_1$};
\node at (78.75:1.6cm) {$\scriptstyle \Phi_1$};
\node at (135:1.075cm) {$\scriptstyle \Phi_2$};
\node at (135:1.7cm) {$\scriptstyle \Phi_2$};
\node at (198:0.9cm) {$\scriptstyle \Phi_1$};
\node at (198:1.6cm) {$\scriptstyle \Phi_1$};
\node at (258.75:0.9cm) {$\scriptstyle \Phi_2$};
\node at (258.75:1.6cm) {$\scriptstyle \Phi_2$};
\node at (315:1cm) {$\scriptstyle \Phi_1$};
\node at (315:1.75cm) {$\scriptstyle \Phi_1$};
\node at (8:0.95cm) {$\scriptstyle \Phi_2$};
\node at (8:1.6cm) {$\scriptstyle \Phi_2$};
\end{tikzpicture}
$$
where the algebras $\Lambda$, $\Lambda_{1}$, $\Lambda_{2}$, $\Lambda_{12}$, $\Lambda_{21}$, $\Lambda_{121}$ are the six NCCRs corresponding to the six crepant resolutions of $\Spec (R)$.
 As explained in \cite[Remark~4.7]{Hirano-Wemyss}, there is a functorial isomorphism $\Phi_1\Phi_2\Phi_1\cong\Phi_2\Phi_1\Phi_2$.

 There are simple modules $\scrS_1$ and $\scrS_2$ for $\Lambda$ with the property that the $\A_{\infty}$-endomorphism algebra $\mathrm{Ext}_{\Lambda}^*(\scrS_1 \oplus \scrS_2)$ agrees with the Floer $A_{\infty}$-algebra of the two core 3-spheres considered previously, cf. \cite[Proposition~3.11]{Smith-Wemyss}.
 By slight abuse of notation, we write $\scrS_1$, $\scrS_2$ for the simple modules for the NCCRs, regardless of which NCCR is under consideration. The labelling is consistent, so that mutation $\Phi_i$ always shifts $\scrS_i$.

Write $C_+$ for the top right chamber in the hyperplane arrangement.

\begin{lem} \label{lem:path bijection}
There is a bijection between morphisms in the Deligne groupoid ending at $C_+$ and bounded $t$-structures on~$\scrC$.
\end{lem}

\begin{proof} This is the analogue of \cite[Corollary 5.3]{August-Wemyss} but dealing with bounded $t$-structures in place of tilting complexes. One considers the map $\alpha \mapsto \Phi_{\alpha}(\bT)$, with $\bT$ the standard $t$-structure, as a~map from morphisms in the groupoid to $t$-structures. Comparing to the proof in op.~cit.\ the injectivity argument goes through exactly as written (using that the pure braid group action on $\scrC$ is faithful, by, e.g., \cite[Corollary~5.5]{August-Wemyss}), whilst the surjectivity argument follows from \cite[Theorem~4.6]{Hara-Wemyss}.
\end{proof}

\begin{lem}\label{lem:same on simples}
If $f$ is a compactly supported symplectomorphism of~$X$ inducing an autoequivalence~$F$, there is a path $\alpha$ so that $F(S_i) = \Phi_{\alpha}(S_i)$ for all the simple objects in the standard heart.
\end{lem}

\begin{proof} We know that $F$ applied to the standard $t$-structure $\bT$ with heart $\heartsuit$ yields some bounded $t$-structure $F(\bT)$ with heart $F(\heartsuit)$, so by the previously established bijection there is some $\alpha$ so that $F(\heartsuit) = \Phi_{\alpha}(\heartsuit)$. Now $F \circ \Phi_{\alpha}^{-1}$ fixes $\heartsuit$ and hence permutes the simples. By Lemma~\ref{lem:trivial on K-theory}, $F$ acts trivially on $K$-theory so induces the identity permutation. On the other hand, the combinatorics of \cite[Proposition 4.8]{Hirano-Wemyss} shows that $\Phi_{\alpha}$ can only act on $K$-theory via one of the six possible matrices realised by the compositions of the standard reflection functors, and one checks directly that only the identity map acts by a permutation matrix. The result follows. \end{proof}

\begin{lem}
If $\Phi_i^{\pm}$ denote the left and right tilts at the simple $S_i$, then $F \circ \Phi_i^{\pm}$ and $\Phi_{\alpha} \circ \Phi_i^{\pm}$ send the standard heart to the same heart.
\end{lem}

\begin{proof}
The definition of the right and left tilts depends only on the heart and choice of simple. The result then follows from Lemma~\ref{lem:same on simples}, which implies that $F \circ \Phi_i^{\pm}$ and $\Phi_{\alpha} \circ \Phi_i^{\pm}$ act in the same way on the torsion pair defining the tilt.
\end{proof}

\begin{lem}
All bounded hearts in $\scrC$ are obtained from a sequence of tilts starting from the standard heart $($i.e.,~the heart exchange graph is connected$)$.
\end{lem}

\begin{proof}
By Lemma~\ref{lem:path bijection}, every bounded t-structure is of the form $\Phi_{\alpha}(\bT)$, where $\Phi_{\alpha}$ is a composition of mutation functors and their inverses. The statement follows using the fact that the mutation functors (and their inverses) give the left and right tilts at simples (see, e.g., \cite[Lemma~5.5]{Hirano-Wemyss}).
\end{proof}

We now induct over the number of tilts needed to reach some given heart from the standard heart. This is a `careless' induction, in the sense that we are not appealing to any normal form for such a sequence of tilts: we are just exhausting all hearts by arbitrary sequences of tilts starting from the standard heart. By this induction, given the $\alpha$ from Lemma~\ref{lem:same on simples} we conclude that $F$ and $\Phi_{\alpha}$ agree on the $t$-structure obtained from $\heartsuit$ by applying any element of the braid group, i.e.,
\[
F \Phi_{\gamma}(\heartsuit) = \Phi_{\alpha}(\Phi_{\gamma}(\heartsuit))
\]
for every path $\gamma$ in the Deligne groupoid ending in the chamber $C_+$. Using the bijection from Lemma~\ref{lem:path bijection}, we see that $F$ and $\Phi_{\alpha}$ agree on all bounded hearts. The following now completes the proof of Theorem~\ref{thm:double bubble case}.

 \begin{cor}
 The association $f \mapsto F$ of Lemma~{\rm \ref{lem:same on simples}} defines a homomorphism
 \[
 \pi_0\Symp_{\rm ct}(X) \to {\rm PBr}_3, \qquad f \mapsto \alpha, \qquad \mathrm{where} \quad F(\heartsuit) = \Phi_{\alpha}(\heartsuit).
 \]
 \end{cor}

\begin{proof}
Given symplectomorphisms $f$, $g$ yielding autoequivalences $F$, $G$ which act trivially on $K$-theory, we define pure braids $\alpha$ and $\beta$ by the property that
\[
F(\heartsuit) = \Phi_{\alpha}(\heartsuit), \qquad G(\heartsuit) = \Phi_{\beta}(\heartsuit).
\]
So then $FG(\heartsuit) = F(\Phi_{\beta}(\heartsuit))$ and we wish to show this equals $\Phi_{\alpha}\Phi_{\beta}(\heartsuit)$. But this follows immediately from the discussion above, which established that $F$ and $\Phi_\alpha$ act the same way on all bounded hearts.
\end{proof}

The categorical/stability condition aspects of the above discussion were not particularly special to the double bubble, but apply more generally to a wide class of 3-fold flopping contractions (the fact that the map $\Phi_{\alpha}$ could not act by an interesting permutation matrix on $K$-theory also applies when there are more simples, by, e.g., the analogue of \cite[Proposition~2.12]{August-Wemyss}).

 As a concrete example, one can consider a `multi-bubble' plumbing $W_f$, which is the Weinstein total space of a Morse--Bott--Lefschetz fibration with fibre $(\bC^*)^2 = T^*T^2$, and with $n+1$ singular fibres having vanishing cycles
 \[
 k_i m + l_i \ell, \qquad k_i = 1 \ \mathrm{or} \ l_i = 1; \qquad (k_i, \pm l_i) \neq (k_j, \pm l_j) \ \mathrm{for} \ 0\leq i\neq j \leq n,
 \]
 where $m$, $\ell$ are a meridian and a longitude in the torus fibre. These spaces were studied by Xie and Li in \cite{Xie-Li}; the hypotheses on the vanishing cycles imply that (i) the core Lagrangians are Lens spaces, and (ii) all the possible surgeries amongst those cores are homology Lens spaces, and not copies of $S^1\times S^2$.
 If
 \[
 f(x,y) = \prod f_j(x,y), \qquad f_j(x,y) = x^{l_j} \pm y^{k_j}
 \]
 and $R_f = \bK[u,v,x,y] / \langle uv-f(x,y)\rangle$, then $\mathrm{Spec}(R_f)$ has an isolated Gorenstein $cA_n$-singularity, which admits a crepant resolution by a chain of $n$ small resolution curves, each with normal bundle $\mathcal{O}_{\bP^1}(-1,-1)$ or $\mathcal{O}_{\bP^1}(0,-2)$.

Let $D_n \subset \bC$ denote the set of $n+1$ critical values of the MBL fibration. There is a relative wrapped category for the pair $(\bC,D_n)$, whose objects are Lagrangians in $\bC\backslash D_n$, where morphism groups are defined using wrapping at the punctures, but where holomorphic curves may cross the divisor $D_n$ and are then counted with weights, yielding a category defined over $R = \bK[t_0,\dots,t_n]$.

 In part building on ideas from \cite{LS}, Xie and Li establish equivalences over any field $\bK$
 \[
\mathrm{Perf}\, \scrW(W_f;\bK) \simeq \operatorname{Perf} \scrW(\bC,D_n) \otimes_R \bK[x,y] \simeq D(Y_f) / \langle \mathcal{O}_{Y_f}\rangle,
 \]
where $R \to \bK[x,y]$ sends $t_i \mapsto f_i(x,y)$ and where the $\bK[x,y]$-linear structure on $\scrW(W_f)$ arises from the map $\bK[x,y] \to SH^0(W_f)^{\times}$ sending $x$, $y$ to the Seidel elements of the fibrewise Hamiltonian circle actions on $W_f \to \bC$. (The fact that this equivalence holds before any completion is a little delicate, and carefully addressed in~\cite{Xie-Li}.)

Exactly as in \cite{Smith-Wemyss}, the relative singularity category is Koszul dual to the `small subcategory' $\scrC_f \subset D(Y_f)$ of objects supported on the exceptional locus and with trivial derived push-forward
\[
D(Y_f) / \langle \mathcal{O}_{Y_f}\rangle \simeq \scrC_f^{!}.
\]
One thus obtains maps
\[
\Auteq_{\bK}(\scrW(W_f;\bK)) \to \Auteq_{\bK}(\scrC_f).
\]
One can then use the classification of hearts, and of objects with no negative $\Ext$'s, from~\cite{Hara-Wemyss} to see that the symplectic monodromy map
\[
{\rm PBr}_{n+1} \to \pi_0\Symp^{\rm gr}(W_f)
\]
splits, just as for double bubbles.

\begin{rmk}
One feature that is special to the double bubble plumbing $X$ (in which the cores and surgeries are actually $S^3$'s with no lens spaces around) is that the action of the pure braid group by symplectomorphisms is known to lift to the compactly supported symplectic mapping class group $\pi_0\Symp_{\rm ct}(X)$; in more complicated geometries (e.g., longer linear plumbing chains) even when the spherical twist generators are compactly supported, it is less clear that the whole representation lifts. Double bubbles have Morse--Bott--Lefschetz fibrations whose singular fibres have non-compact critical set, but at least naively the braid relations in that viewpoint hold at the level of maps which may have non-proper support. The fact that the braid relation holds at the level of compactly supported maps for $X$ itself is proved in~\cite{Smith-Wemyss} using our alternative presentation of~$X$ as an affine flag $3$-fold.
\end{rmk}

 \begin{rmk}\label{rmk:multibubble as surface category}
 It would be interesting to directly infer splitting from the identification of $\scrW(W_f)$ with the category $\scrW(\bC,D_n) \otimes_R \bK[x,y]$, which gives a trivial proof of \emph{faithfulness} of the pure braid action.
 \end{rmk}

In all these cases, the underlying monodromy groups are classical (pure) braid groups. We know that more interesting, e.g., non-Coxeter, groups arise from more interesting flopping contractions, though we do not always have candidate geometric mirrors.

\begin{qn}
Let $X$ be mirror to the small resolution of a $3$-fold flopping contraction $Y \to \mathrm{Spec}(R)$, and let $\mathcal{M}_X$ denote the complexified hyperplane complement with fundamental groupoid the pure braid group ${\rm PBr}(\scrC)$ associated to the small category $\scrC$. Does symplectic monodromy split over $\pi_1\mathcal{M}_X$?
\end{qn}

\section{The conifold and its inexact deformations} \label{sec:conifold}

The double bubble has a cousin, given by the `conifold quiver'
\[
\xymatrix@C=1em{\bullet \ar@/^0.5pc/@{=>}[rr]^{e,e'} && \bullet, \ar@/_-0.5pc/@{=>}[ll]^{f,f'} } \qquad W = efe'f' - ef'e'f
\]
with two vertices and two arrows in each direction. On the $B$-side this again is well-known as describing the derived category of sheaves on the small resolution of the ordinary double point, supported on the exceptional curve $\bP^1$. A better mirror geometry is provided by the space $Y = (\mathcal{O}(-1) \oplus \mathcal{O}(-1)) \setminus C$, where $\mathcal{O}(-1) \oplus \mathcal{O}(-1)$ is the small resolution of the threefold ordinary double point and $C$ is a smooth quadric in the complement of the exceptional curve. The space $Y$ still admits a morphism $f\colon Y \to \bP^1$, given by the restriction to $Y$ of the vector bundle projection. To $Y$ one associates two categories
\[
\scrC = \{E^{\bullet} \in D(Y) \mid Rf_*(E^{\bullet}) = 0\} \qquad \mathrm{and} \qquad \scrD = \{E \in D(Y) \mid \mathrm{supp}\,H^*(E^{\bullet}) \subset \bP^1\}
\]
where $D(Y)$ denotes the usual bounded derived category of coherent sheaves. For a symplectic model, following \cite{CPU} and \cite{Keating-Smith}, we consider the affine variety~$X$ which is the complete intersection
\begin{equation} \label{eqn:conifold mirror}
X = \bigl\{ (u_1,v_1,u_2,v_2,z) \in \bC^2\times \bC^2 \times \bC^* \mid u_1v_1=z-1, \, u_2v_2=z+1\bigr\}.
\end{equation}
This is also a plumbing\footnote{Where the double bubble was a Morse--Bott--Lefschetz fibration over $\bC$ with general fibre $T^*T^2$ and 3 singular fibres, associated to Bott vanishing circles $a$, $b$, $a+b$ with $a$, $b$ the longitude and meridian for the torus, cf.\ Figure~\ref{Fig:vcycles}, the conifold mirror is a Morse--Bott--Lefschetz fibration over $\bC^*$ with the same general fibre, and with Bott vanishing circles $a$, $b$ over the two critical fibres. The two core 3-spheres are now associated to paths between the two critical values either side of the puncture at the origin in $\bC^*$.} of two $T^*S^3$'s along a Hopf link, see \cite[Lemma~3.3]{Keating-Smith}.
Write $\scrF_{\mathrm{core}}(X)$ for the subcategory of the compact Fukaya category split-generated by the core Lagrangian spheres of the plumbing. Mirror symmetry now states the following.
\begin{thm}[{\cite[Theorems 1.2 and 1.4]{CPU}}]
There are equivalences $\scrF_{\mathrm{core}}(X) \!\simeq \! \scrC$ and ${\scrW(X) \!\simeq\! \scrD}$.
\end{thm}

The manifold $X$ admits a MBL fibration over $\bC^*$, by mapping to the $z$ coordinate, with fibre $\bC^*\times \bC^*$; singular fibres at $\pm1$; and core matching spheres associated (up to isotopy) to the matching paths given by the upper and lower half-circles of the unit circle in $\bC^\ast$. Varying the fibration gives a family $\scrC \to \mathrm{Conf}_2(\bC^*)$ and hence a representation
\[
{\rm PBr}_3 \to \pi_0\Symp(X).
\]
This no longer has image in $\Symp_{\rm ct}(X)$. Forgetting the first respectively third strands defines a natural map
\[
{\rm PBr}_3 \to {\rm PBr}_2 \times {\rm PBr}_2
\]
whose kernel we denote ${\rm PBr}_3^c$; these maps act by compactly supported symplectomorphisms. One has that ${\rm PBr}_3^c = \bZ^{\ast \infty}$, a countably infinite rank free group. For details, see \cite[Section 3]{Keating-Smith}.

\begin{thm}[Keating and Smith]
The monodromy homomorphism
\[
\rho\colon \ {\rm PBr}_3^c \to \pi_0\Symp_{\rm ct}(X)
\]
is split.
\end{thm}

\begin{proof} This is spelled out explicitly in \cite[Theorem 1.1]{Keating-Smith}. \end{proof}

This example illustrates that splitting can have non-trivial implications, e.g.,~here it yields a~Stein manifold with infinitely generated symplectic mapping class group.

\begin{figure}
 \centering
 \includegraphics[width=0.7\linewidth]{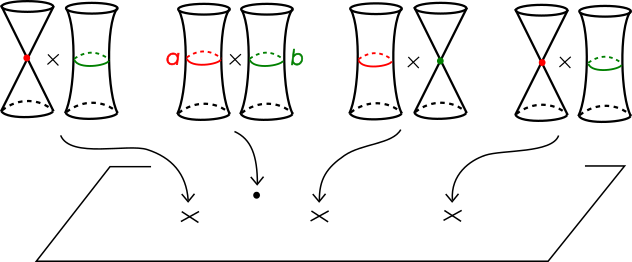}
 \caption{Morse--Bott--Lefschetz fibration on $(W_0, \omega_0)$. The crosses in the base denote singular values, and the dot a regular one.
 \label{fig:MLB-for-W_0}}
\end{figure}

There is a geometry which mediates between the double bubble and the conifold. Let $(W_0,\omega_0)$ be the Stein manifold which is a plumbing of two copies of $T^*S^3$ along a circle so that the Lagrange surgery of the core spheres is an $S^1\times S^2$. This is the Morse--Bott--Lefschetz fibration with fibre $(\bC^*)^2$ and, in the notation from Figure~\ref{Fig:vcycles}, vanishing cycles $a,b,a$ respectively. See Figure~\ref{fig:MLB-for-W_0}.

\begin{rmk}
 The double bubble is a plumbing of two 3-spheres along a circle so that the Lagrange surgery is again a 3-sphere. The space $W_0$ is a plumbing of two copies of $T^*S^3$ along a circle so that the Lagrange surgery is $S^1\times S^2$. The conifold mirror $X$ is a plumbing of two copies of $T^*S^3$ along a Hopf link; its universal cover is a plumbing of countably many copies of $T^*S^3$ along circles in a linear chain, such that the Lagrange surgery of any adjacent pair yields an $S^1\times S^2$. This is the sense in which $W_0$ `mediates between' the double bubble and conifold mirror, and why inexact deformations of $W_0$ provide a template for studying inexact deformations of the conifold mirror.
\end{rmk}

\begin{rmk}\label{rmk:classify 3-spheres}
 The analysis from \cite{Smith-Wemyss} did not yield a classification of Lagrangian 3-spheres in~$W_0$, in contrast to the case for the closely related double bubble plumbing from Section~\ref{sec:double-bubble} (called~$W_1$ in op.~cit.). The space~$W_0$ is an affine quartic, and the total space of a conic fibration over~$T^*S^2$, which in turn is a conic fibration over a plane. Note that~$W_0$ admits a $\bC^*$-action, and any Lagrangian sphere defines an equivariant module, cf.~\cite[Proposition~6.5]{Auroux-Smith} which follows \cite[Corollary~A.12]{Seidel-dilating}. Recent speculations of Lekili and Segal~\cite{LS} give a description of the $(\bC^*)^2$-equivariant Fukaya category of~$W_0$, for the fibrewise action, in terms of a partially wrapped Fukaya category of a thrice-punctured plane. Starting from this, and using classification results for objects in Fukaya categories of surfaces, one could hope to classify Lagrangian 3-spheres in~$W_0$. This would be a rather different route to classification than that taken for $W_1$ in \cite{Hara-Wemyss, Smith-Wemyss}, closer to the recent arguments of~\cite{Xie-Li} (cf.~Remark~\ref{rmk:multibubble as surface category} above).
\end{rmk}

One easily checks that $H^2(W_0;\bZ) \cong \bZ$, so there is an inexact deformation $\omega_t$ of the symplectic form $\omega_0$ on $W_0$. In fact, this can be realised geometrically. Note that the original Morse--Bott--Lefschetz fibration on $W_0$ has monodromies on $\bC^*\times \bC^*$ which are given by
\[
\tau_a \times \{{\rm id}\}, \quad \{{\rm id}\} \times \tau_b, \quad \tau_a \times \{{\rm id}\},
\]
respectively, where $\tau_{\bullet}$ denotes the Dehn twist in the curve $\bullet$. The corresponding vanishing cycles, for appropriate vanishing paths, are $\{a,b,a\}$ in a hopefully obvious shorthand. Let $\phi_x$ denote the symplectomorphism of $\bC^* = \bR\times S^1$ which is given by a translation by $x$ in the $\bR$-factor, i.e.,~the symplectomorphism (unique up to Hamiltonian isotopy) with flux $x \in H^1(\bC^*;\bR)\cong \bR$.

\begin{lem}
The space $(W_0,\omega_t)$ admits a Morse--Bott--Lefschetz fibration with the same monodromies as before, but the vanishing cycles associated to the previously chosen vanishing paths are now $\{a,b,(\phi_t\times \id)(a)\}$.
\end{lem}

\begin{proof}
 Cut $W_0$ along the codimension one hypersurface which is the preimage of a vertical line in $\bC$ disjoint from the critical values of the MBL fibration, and such that the two critical values with vanishing cycle $a$ lie on either side of the cut; and reglue by the map given by translation on the first factor of the fibre and the identity on the second. This is a symplectomorphism, so the result is again a symplectic manifold, manifestly diffeomorphic to $W_0$ (we have only changed the locally Hamiltonian structure of the fibration; the re-gluing is by a map smoothly isotopic to the identity). Since the change in symplectic structure is localised away from the critical submanifolds, the local monodromies are unchanged, but at least one of the vanishing paths will cross the cut locus and the corresponding vanishing cycle will therefore change as indicated.
\end{proof}

\begin{cor}\label{cor:disjoinable}
 Let $0 \neq t \in \bR$. In the inexact deformation $(W_0,\omega_t)$, the two core Lagrangian $3$-spheres are disjoint.
\end{cor}

\begin{proof}
 The two matching spheres now meet the central singular fibre in copies of a circle $S^1 \times \{pt\} \subset \bC^*\times (\bC\vee\bC)$, where the (first factor) circles for the two different matching spheres differ by a translation of flux $t$.
\end{proof}

Note that the Lagrangian $S^1 \times S^2 \subset (W_0,\omega_0)$ which was fibred over a matching path between the outermost critical values is no longer Lagrangian, since the two-sphere factor has acquired positive area. The space $W_0$ carries an action by symplectic parallel transport of the mixed braid group $M{\rm Br}_3$. Consider the generator given by the lift of a half-twist, in $\bC$, along a matching path between the two critical values of type $a$. Observe that this is also naturally a lift under a~map $W_0 \to T^*S^2$ (forgetting the second coordinate on each fibre) of the Dehn twist in the $S^2$, which no longer has infinite order after the inexact deformation of $\omega$. Similarly, the twists in the core spheres now commute, since these are now disjoint.

The conifold mirror $X$ introduced in \eqref{eqn:conifold mirror} has $H^2(X;\bR) \cong \bR^2_{s,t}$. Let $(X, \omega_{s,t})$ be a corresponding inexact deformation of the Weinstein manifold $X$. We obtain a geometric model for $(X,\omega_{s,t})$ by considering the Morse--Bott--Lefschetz fibration over $\bC^*$ with two singular fibres with monodromies
\[
\tau_a \times \id, \ \{\phi_s\} \times \{\phi_t\} \ \mathrm{and} \ \id \times \tau_b
\]
around $-1$, $0$ and $ 1$, respectively, for the obvious three loops based at $i \in \bC^\ast$.
Note that, as one iterates monodromies around the critical fibres, the fluxes add, so one finds the following.

\begin{cor}
 Let $s$, $t$ both be non-zero. The deformation $(X,\omega_{s,t})$ of the conifold mirror contains infinitely many pairwise disjoint Lagrangian 3-spheres.
\end{cor}

\begin{figure}
 \centering
 \includegraphics[width=0.7\linewidth]{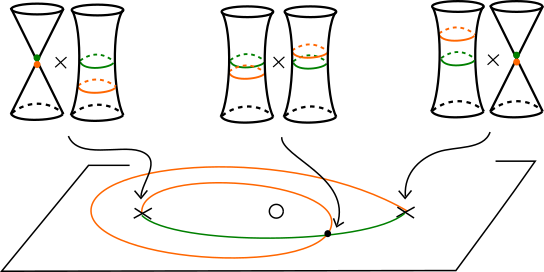}
 \caption{Example of disjoint Lagrangian matching spheres in $(X, \omega_{s,t})$.}
 \label{fig:Inexact-conifold-with-spheres}
\end{figure}

\begin{proof}
 Consider matching paths between the two critical values with larger and larger winding number around the origin; given two distinct such, at any intersection point of the base matching paths, the winding number difference means at least one of the fibrewise factors has been translated by a non-trivial flux. See Figure~\ref{fig:Inexact-conifold-with-spheres}.
\end{proof}

\begin{rmk}
 The Lagrangian spheres associated to matching paths are still the only Lagrangian spheres we can `see' after deformation, but we no longer have even a tentative route to a~classification guaranteeing that these are the only spherical objects.
\end{rmk}

Each Lagrangian matching sphere has an associated Dehn twist. Disjointness implies these twists now commute. By considering the action on the set of graded Lagrangian spheres, and using that $\tau_L(L)=L[-1]$ but $\tau_L$ does not shift $K$ if $K \cap L = \emptyset$, one sees that the homomorphism from $\bZ^{\oplus \infty} \to \pi_0\Symp_c^{\rm gr}(X,\omega_{s,t})$ will be injective when $s,t \neq 0$.
\begin{qn}
 Does the resulting homomorphism $\bZ^{\oplus \infty} \to \pi_0\Symp_c(X,\omega_{s,t})$ split?
\end{qn}

\begin{rmk}
 The approach suggested in Remark \ref{rmk:classify 3-spheres} does not obviously survive the deformation of the symplectic form, since this has no obvious interpretation in the model for the equivariant category formulated on the reduced space of the $\bC^*$-actions.
\end{rmk}

 As a digression, we note that one can `compactify' the preceding example:

 \begin{prop}
 There is a compact symplectic six-manifold which contains infinitely many pairwise disjoint Lagrangian $3$-spheres.
 \end{prop}

\begin{proof}
 We consider a Morse--Bott--Lefschetz (MBL) fibration with fibre $T^2\times T^2$ constructed as follows. Begin with the four-manifold which is a symplectic $T^2$-fibration over a torus with monodromies $\id$ and $\phi_s$ for some flux $s$. This is locally Hamiltonian, so admits a symplectic structure. We fibre sum this with a rational elliptic surface, to obtain a symplectic fibration $Z \to T^2$ over the two-torus with a dozen nodal singular fibres lying over a disc $D \subset T^2$, still with monodromies $\{id,\phi_s\}$.

 We now take two copies $Z$ and $Z'$ of this four-manifold, where $Z$ has monodromies $\{\id,\phi_s\}$ and $Z'$ has monodromies $\{\phi_t,\id\}$ (so the flux has switched to the other factor of the base). We also insist that the singular fibres for $Z$ (resp.~$Z'$) lie over disjoint discs $D,D' \subset T^2$. Now take the fibre product $Z \times_{T^2} Z'$. The disjointness of $D$ and $D'$ means that at most one of the two input fibrations has a singular point over a given $p\in T^2$, which ensures that the fibre product is the total space of a MBL fibration with general fibre~$T^4$.

 We now impose the condition that the non-trivial fluxes $s,t \in H^1(T^2;\bR) \cong \bR^2$ are irrational, so have infinite order in the symplectomorphism group of the fibre. One can construct Lagrangian matching spheres in the total space exactly as in the inexact conifold. The irrationality of the fluxes ensures that the spheres associated to matching paths which meet in the base, as in Figure~\ref{fig:Inexact-conifold-with-spheres}, are still disjoint in the total space. This is because the iterated translations of a compact two-torus in $T^4$ by an irrational angle in a complementary direction (i.e., writing $T^4 = T^2 \times T^2$ and translating the first factor by iterations of an irrational flux in the second factor) remain pairwise disjoint for all iterations.
\end{proof}

Such examples do not seem to have been observed previously.

\section{The birational projective plane}\label{sec:birational-proj-plane}

 Let $X$ denote the `universal mirror' to log CY surfaces constructed in \cite[Section 2.5]{Keating-Ward} (where it is denoted $M_\mathrm{univ}$). It is a Stein manifold of infinite type, obtained by attaching an infinite sequence of 2-handles to $T^*S^2$, one attached along the Legendrian conormal for every isotopy class of oriented simple closed curve on $T^2$. It is also the total space of a Lagrangian torus fibration over $\bR^2$ with countably many singular fibres, over $\bZ^2\backslash\{0\}$, with invariant directions all passing through the origin (this ensures that the symplectic form is exact). In abuse of notation, we will write
 \[
 \pi_0\Symp_{\rm ex}(X) = \Symp_{\rm ex}(X)/\Ham_{\rm ct}(X)
 \]
 for the quotient of \emph{exact} (but not necessarily compactly supported) symplectomorphisms by compactly supported Hamiltonian isotopy. Let $\Omega = {\rm d}x/x \wedge {\rm d}y/y$ denote the standard meromorphic volume form on $\bP^2$. The volume-preserving Cremona group $\Bir\bigl(\bP^2,\Omega\bigr)$ is the group of birational maps which preserve~$\Omega$. Let $\Bir_e\bigl(\bP^2,\Omega\bigr)$ be the subgroup of the volume-preserving Cremona group which is generated by ${\rm SL}_2(\bZ)$ together with the basic birational cluster transformation~$E$, with
 \[
 E\colon \ (x,y) \mapsto \bigl(x,y(1+x)^{-1}\bigr).
 \]
 (The whole $\Bir\bigl(\bP^2,\Omega\bigr)$ is generated by this subgroup and $(\bC^\ast)^2$, acting by multiplication coor\-dinate-wise.)
\cite[Section~3.2]{Keating-Ward} constructs a natural map $\Bir_e\bigl(\bP^2,\Omega\bigr) \to \pi_0\Symp_{\rm ex}(X)$.

\begin{rmk}
 For a field $k$, the group $\Bir\bigl(\bP^2(k),\Omega\bigr)$ certainly depends on $k$, for instance it contains the multiplicative torus $(k^*)^2$. However, as Blanc points out in~\cite{Blanc}, the group $\Bir_e\bigl(\bP^2,\Omega\bigr) = \langle {\rm SL}_2(\bZ), E\rangle$ (which is independent of $k$) embeds in $\Bir\bigl(\bP^2(k),\Omega\bigr)$ for any $k$.
\end{rmk}

As noted above, the space $X$ is the total space of a Lagrangian torus fibration over $\bR^2$ with infinitely many singular fibres. Let $\Symp_{\rm fib}(X) \subset \Symp_{\rm ex}(X)$ denote the subgroup of maps which preserve the set of Hamiltonian isotopy classes of smooth Lagrangian torus fibres near infinity, i.e.,~for which there is an open neighbourhood of infinity (depending on the particular map $\phi$) such that $\phi$ takes fibres in that neighbourhood to tori which are Hamiltonian isotopic to fibres. Viewing the fibres near infinity as the pre-Lagrangian tori which live inside the mapping torus boundary of a large circle of regular values, this can be seen as a version of being `contact type' at infinity in this special case; compactly-supported maps, linear maps and the cluster transformation $E$ naturally lift to $\Symp_{\rm fib}(X)$.

 \begin{conj}
 \label{conj:birational plane case}
The map $\Bir_e(\bP^2,\Omega) \to \pi_0\Symp_{\rm fib}(X)$ splits.
 \end{conj}

 We will outline an argument for this, but a full proof would need more care because we use inexact Lagrangians (general Lagrangian torus fibres of an almost toric fibration), including some properties expected by experts -- whereas the set-up and HMS proofs of \cite{Keating-Ward} restrict to the exact Fukaya category.
 We have packaged these properties into Assumption~\ref{assump:inexact}, and sketch a proof of Conjecture~\ref{conj:birational plane case} subject to the assumption.

 The broad approach to Conjecture~\ref{conj:birational plane case} is to associate toric charts $(\Lambda_{\bC}^*)^2$ in the mirror space $U$ to $X$ to a particular Lagrangian torus $T$ and its symplectomorphism image $\phi(T)$, and to see that these charts differ by a birational map, which is the one we associate to $\phi$ (the precise argument differs from this somewhat). We begin by describing $X$ and its mirror space $U$ as direct limits of finite type spaces.

 Since the space $X$ admits a handle decomposition in which countably many two-handles are attached to $T^*T^2$, the zero-section of that $T^*T^2$ defines an exact Lagrangian torus $T\subset X$.
 The space $X$ is constructed as the direct limit of finite-type subdomains $X_{\xi}$, where the $\xi$ run over a~directed system $\Xi$.
 This indexes log Calabi--Yau surfaces $(Y_\xi, D_\xi)$ with maximal boundary and split mixed Hodge structure, equipped with explicit toric models (see \cite[Definition 2.4]{Keating-Ward}). We have
\begin{equation*}
 X = \varinjlim_{\xi \in \Xi} X_{\xi}
\end{equation*}
and the mirror is a corresponding limit of open log CY surfaces:
\begin{equation*}
U = \varinjlim_{\xi \in \Xi} U_{\xi},
\end{equation*}
where $U_{\xi} = Y_{\xi} \setminus D_\xi$.

 \begin{lem} \label{lem:torus-mirror-over-C}
 For any $\xi \in \Xi$, there is a quasi-isomorphism of $A_\infty$-categories
\begin{equation*}
 \Upsilon\colon \ \ccoh(U_{\xi};\bC) \simeq \tw \scrW(X_{\xi};\bC),
\end{equation*}
$\ccoh (U_\xi; \bC)$ denotes a dg enhancement of the derived category of coherent sheaves $D^b \Coh (U_\xi; \bC)$ $($recall these are unique up to quasi-isomorphism {\rm \cite{Lunts-Orlov})}, and $\tw \scrW(X_{\xi};\bC)$ is the category of twisted complexes enlarging $\scrW(X_{\xi};\bC)$.
 Moreover, for any choice of brane structure $\rho$ on~$T$,
 $\Upsilon$ maps $(T, \rho)$ to the skyscraper sheaf of a closed point in $U_\xi$ (up to quasi-isomorphism).
 \end{lem}

\begin{proof} See the equivalence of categories \cite[Section 2.3]{Keating-Ward}, which in turn imports this from \cite{Hacking-Keating-1}.
The claim about the mirror to $(T, \rho)$ is in \cite[Lemma 2.24]{Hacking-Keating-3}.
\end{proof}

For what follows, it will be important to pass from the exact Fukaya category over $\bC$ to the more general tautologically unobstructed Fukaya category, which is defined over the one-variable Novikov field $\Lambda_{\bC}$. Objects of the tautologically unobstructed category include pairs $(L,J_L)$, where $L$ is a compact Maslov zero Lagrangian, and $J_L$ is a choice of compatible almost complex structure for which $L$ bounds no $J_L$-holomorphic discs; this is a generic property in the almost complex structure, for virtual dimension reasons, in particular every Maslov zero $L$ admits such $J_L$.
We remind the reader that $\Lambda_{\bC}$ is an algebraically closed field of characteristic zero.

The log Calabi--Yau surfaces $U_\xi$ are characterised by a split mixed Hodge structure condition which is naturally formulated over $\bC$.
\begin{lem}
 The log CY pairs $(Y_\xi, D_\xi)$, $U_\xi$, and the direct system of inclusions of the $U_\xi$ indexed by $\xi \in \Xi$, are all naturally defined over $\Lambda_{\bC}$.
\end{lem}

\begin{proof} Base-change for the inclusion $\bC \to \Lambda_{\bC}$.
\end{proof}

\begin{assum}\label{assump:inexact} There is an equivalence
 \[
 \psi\colon \ \tw \scrW(X_{\xi};\Lambda_{\bC
}) \stackrel{\sim}{\longrightarrow} \mathrm{coh} (U_{\xi};\Lambda_\bC)
 \]
 with the following properties:
\begin{enumerate}\itemsep=0pt
 \item[$(i)$] Under $\psi$, the local systems on $T$ define a chart $(\Lambda_{\bC}^*)^2 \hookrightarrow U_\xi$.
 \item[$(ii)$] For any sufficiently large compact subset $Z \subset \bR^2$, under $\psi$ each smooth Lagrangian torus fibre of $X_{\xi} \to \bR^2$ lying over $\bR^2 \setminus Z$ is taken to the skyscraper sheaf of a point in $U_{\xi}$. Furthermore, such points cover a Zariski open subset of a punctured neighbourhood of the boundary $D_{\xi} = Y_{\xi} \setminus U_{\xi}$.
\end{enumerate}
\end{assum}

{\bf Discussion.}
For (i), note that by~\cite{GPS}, we can pick a finite set of split-generators for $\scrW(X_{\xi})$.
By expressing the exact torus $T$ in terms of those split-generators, $T$ gives rise to an algebraic family of modules over the Fukaya category parametrised by $H^1(T;\bG_m)$; see \cite{Seidel:Flux} and \cite[Sections~2 and~3]{Kartal:algebraic}. The conclusion of the part (i) is then essentially \cite[Note~1.12]{Kartal:algebraic} (though this is only sketched).
Here is a brief outline: the mirror equivalence takes the algebraic family of modules over $\tw \scrW(U_{\xi})$ to an algebraic map from $H^1(T;\bG_m)$ to the Artin stack parametrising the moduli of objects of $\ccoh(U_\xi)$. In general the moduli of objects in such a $dg$-category is only locally of finite type \cite{TV}, and can be quite pathological. The content of \cite{Kartal:algebraic} is to nonetheless identify non-unitary deformations with flux images, by crucially using the group-like property of the family. 

\emph{The rest of this section is contingent on Assumption~{\rm \ref{assump:inexact}}}, but we will not spell this out in every statement.

 We wish to consider the corresponding charts for tori $T$ and $\phi(T)$, where $\phi \in \pi_0\Symp_{\rm ex}(X)$.

 \begin{lem}
 Let $T' \subset X$ be any exact Lagrangian torus. Then there is some $\xi = \xi(T')$ such that $T' \subset X_{\xi} \subset X$. \end{lem}

 \begin{proof} The paper \cite{Keating-Ward} constructs an almost-toric torus fibration $X \to \bR^2$ (called $\pi_{\textrm{univ}}$), which~is a limit of almost toric fibrations from the $X_{\xi}$ to increasingly large subsets of $\bR^2$ (see \cite[Sections~2.2 and~2.3]{Keating-Ward} for details). As $T'$ projects to a compact subset, it must lie inside some~$X_{\xi}$.
 \end{proof}

Exactness implies that each of $T$ and $T'$ support a family of objects of $\scrW(X_{\xi}; \Lambda_{\bC})$ para\-metri\-zed by $(\Lambda_\bC^*)^2$, since their self Floer cohomology is non-trivial for any rank one Novikov-valued local system.
Under the given equivalence $\tw \scrW (X_{\xi} ; \Lambda_{\bC}) \to \ccoh(U_{\xi}; \Lambda_{\bC})$,
the family supported on $T$ defines a chart $(\Lambda_\bC^*)^2 \hookrightarrow U_\xi$, by Assumption \ref{assump:inexact}.
The corresponding family of objects parametrised by $T'$ are in principle only objects of the derived category: they need not \emph{a priori} be realised as another chart.

Recall that $U_{\xi} = Y_{\xi} \setminus D_{\xi}$, where $(Y_\xi, D_\xi)$ is a log CY surface with maximal boundary. That is, $Y_\xi$ is a smooth rational projective surface, and $D_\xi \in |{-K_{Y_\xi}}|$ is either an irreducible nodal rational curve or a cycle of $k\geq 2$ smooth rational curves.
Maximality of the boundary (i.e., that the dual graph to $D_{\xi}$ has a zero-stratum) implies that the components are homologically non-torsion (this is clear if the boundary is irreducible, since a rational surface never has torsion canonical class, and follows by considering the intersection form on components of the boundary when it is a cycle of rational curves), so $K_{Y_{\xi}}$ is non-torsion in Chow. We have the following.

\begin{lem}\label{lem:point-like support}
Fix an algebraically closed characteristic zero field $k$. If $Y$ is a projective surface over $k$ with $K_Y \neq 0 \in CH^*(Y)_{\bQ}$, then a point-like object in $D(Y)$ is either the skyscraper sheaf of a point up to shift, or is supported on a~$($possibly reducible$)$ connected closed curve~$C$ with $K_Y \cdot C = 0 \in \bZ$.
\end{lem}

\begin{proof} This is part of \cite[Lemma~2.2]{Hu-Krah} (note that the proof is written over $\bC$, but works identically over~$k$).
\end{proof}

We are now ready to sketch the proof of Conjecture~\ref{conj:birational plane case}. Suppose we start with a class $[\phi] \in \pi_0\Symp_{\rm fib}(X)$, with representative~$\phi$, a symplectomorphism assumed to preserve fibres at infinity. Set $T' = \phi(T) $.
The key difficulty in establishing any version of Conjecture~\ref{conj:birational plane case} is to rule out the case in which the complexes associated to the family of objects on~$T'$ have support containing a~$(-2)$-curve in~$U_{\xi}$; this will be where fibred-ness enters.

Continuing to work over $\Lambda_\bC$, we have the following.

\begin{lem}
The objects of $D(U_{\xi})$ associated to local systems on $T'$ have compact support contained in $U_{\xi}$.
\end{lem}

\begin{proof}
We know $T'$ has trivial Floer cohomology with all the Lagrangian torus fibres of $X_{\xi} \to \bR^2$ outside some compact subset $Z \subset \bR^2$. The result then follows from Assumption~\ref{assump:inexact}\,(ii).
\end{proof}

\begin{cor}
 The tori $T$ and $T'=\phi(T)$ define two families of point-like objects in $\Coh(U_\xi)$ associated to two charts $(\Lambda_\bC^*)^2 \dashrightarrow U_{\xi}$, which differ by a well-defined birational automorphism $\Psi_{\xi}(\phi)$.
\end{cor}

\begin{proof}[Sketch]
By hypothesis, we begin with a globally defined family $(\Lambda_{\bC}^*)^2 \hookrightarrow U_{\xi}$ of points which are the mirrors to local systems on the preferred torus $T$.

There is a family $\EuI = \{\EuI_t\}_{t \in (\Lambda_\bC^*)^2}$ of complexes of modules over $U_{\xi}$ parametrized by the local systems on $\phi(T) = T'$. The fundamental results of Kartal \cite{Kartal:algebraic} imply that this family is algebraic.

By hypothesis $\phi$ is ``fibred at infinity'', meaning that its mirror automorphism sends skyscrap\-er sheaves of points near infinity to skyscraper sheaves of points. Take such an image point $p \in U_{\xi}$ not lying on one of the countably many $(-2)$-curves. One sees that there is some parameter $t_0$ for which $\scrJ_{t_0}$ has support which does not contain any $(-2)$-curve, so is in fact a point, and is thus quasi-isomorphic to the skyscraper sheaf $\mathcal{O}_p$
(cf.~the discussion of Assumption \ref{assump:inexact}\,(i)).

For any $(-2)$-curve $C \subset U_{\xi}$, and point $q \in C$, we have $\Ext^*(\mathcal{O}_p,\mathcal{O}_q)$ = 0. Semicontinuity of the rank of $\Ext$-groups, and the algebraicity of the family of modules $\EuI$, means that the rank of $\Ext^*(\EuI_t,\mathcal{O}_q)$ jumps up on at most a Zariski-closed subset of parameters $t \in (\Lambda_{\bC}^*)^2$. Thus, for any fixed $C$ and $q\in C$, there is a Zariski open $W_q \subset (\Lambda_{\bC}^*)^2$ with the property that
\begin{equation}\label{eqn:first Ext}
\Ext^*(\EuI_t,\mathcal{O}_q) = 0 \qquad \forall t \in W_q.
\end{equation}

Suppose for contradiction that there is a Zariski open subset $V \subset (\Lambda_{\bC}^*)^2$ of parameters with the property that for $t\in V$, the support $\mathrm{supp}(\EuI_t) \supset C_t$ contains some perhaps reducible $(-2)$-curve $C_t$ (which may depend on the point $t\in V$). Since there are only countably many $(-2)$-curves, and a Zariski open cannot be covered by a countable family of positive codimension subvarieties, there must be at least one curve $C_0$ appearing in the support of a Zariski open set $V_{C_0}$ of the family of modules. Now pick any point $q_0\in C_0$. Then
\begin{equation} \label{eqn:second Ext}
\Ext^*(\EuI_t,\mathcal{O}_{q_0}) \neq 0 \qquad \forall t \in V_{C_0}.
\end{equation}

But now $V_{C_0} \cap W_{q_0} \neq \emptyset$, and equations \eqref{eqn:first Ext}, \eqref{eqn:second Ext} then yield a contradiction.

It follows that the Lagrangian torus $T'$ defines a birational chart $(\Lambda_{\bC}^*)^2 \dashrightarrow U_{\xi}$, i.e.,~there is a~Zariski open subset of points of the domain which are mapped to skyscraper sheaves. Since~$U_{\xi}$ is irreducible, the charts associated to $T$ and $T'$ necessarily differ by a birational map.

This birational map is only well-defined up to the ambiguity coming from a choice of identification of the local systems on each of the tori with $(\Lambda^*)^2$, which amounts to the ambiguity of composition with a linear map. The hypothesis that $T'=\phi(T)$ for $\phi \in \Symp(X)$ gives a distinguished identification between local systems on $T$ and on $T'$. This implies that the birational map $\Psi_{\xi}(\phi)$ between the two charts is uniquely defined.
\end{proof}

\begin{lem}
The birational map $\Psi_{\xi}(\phi)$ preserves $\Omega$.
\end{lem}

\begin{proof} This is a consequence of Poincar\'e duality in Floer cohomology, compare to \cite[Section~11]{Seidel:categorical_dynamics}.
To spell it out somewhat, suppose we have
two exact Lagrangian tori $T, T' \subset X$, and an isomorphism in the
Fukaya category between Zariski open subsets of the families of local systems on $T$ and $T'$. The tangent space to the
family of local systems can be identified at a given point with $HF^1$, and the isomorphism
$HF^*(T,T) \cong HF^*(T',T')$ is one of rings. Poincar\'e duality implies that the isomorphism is compatible with the trace map $\int_{T'}\colon H^2(T',T') = HF^2(T',T') \to \Lambda_\bC$, since if $x \in HF(T,T')$ and $y \in HF(T',T)$
are mutually inverse isomorphisms, then duality says
\begin{equation} \label{eqn:trace}
\int_T xay = \int_{T'} yxa = \int_{T'} a
\end{equation}
for $a\in HF^2(T',T')$ (and similarly with the roles of $T,T'$ reversed, or for the isomorphisms at a pair of non-trivial local systems).
This in turns means that the isomorphism on $HF^1 = H^1$ has determinant one. Going to the $B$-side, $HF^1(T,T)$ is identified with $\Ext^1(\cO_p,\cO_p)$, which is the holomorphic tangent space to $U_{\xi}$ at~$p$. In the co-ordinate chart $p \in (\Lambda_{\bC}^*)^2 \subset U_\xi$ the trace~\eqref{eqn:trace} is the map on pairs of tangent vectors $(v,v') \mapsto \Omega(v,v')$, for the holomorphic volume form $\Omega = {\rm d}x/x \wedge {\rm d}y/y$. The determinant one condition on $HF^1$ thus means that the induced birational map on $U_\xi$ is volume-preserving, as required.
\end{proof}

Up to this point, one can take any `sufficiently large' subdomain $X_{\xi}$ of $X$ or sufficiently large log CY $U_{\xi}$.

\begin{lem}
The birational map $\Psi_{\xi}(\phi)$ depends only on $\phi$, i.e.,~it is unchanged under embedding $U_{\xi} \hookrightarrow U_{\xi'}$.
\end{lem}

\begin{proof} The birational map is, by construction, determined by knowing which pairs of local systems $\zeta \to T$ and $\zeta' \to T'$ define quasi-isomorphic objects in $\scrF(X_\xi)$. The existence of a quasi-isomorphism is witnessed by certain holomorphic curves in $X_{\xi}$, which continue to exist under an enlargement of Stein domains $X_{\xi} \hookrightarrow X_{\xi'}$. On the other hand, a monotonicity argument shows that if $\zeta$ and $\zeta'$ are quasi-isomorphic in $X_{\xi'}$, then by suitably stretching the neck one can confine the holomorphic curves underlying that quasi-isomorphism to $X_{\xi}$, so they were already quasi-isomorphic there.
\end{proof}

We therefore just write $\Psi(\phi)$ for the birational map $\Psi_{\xi}(\phi)$ for any choice of sufficiently large~$\xi$.

\begin{lem}\label{lem:done}
The association $\phi \mapsto \Psi(\phi)$ is a homomorphism $\pi_0\Symp_{\rm fib}(X) \to \Bir_e\bigl(\bP^2,\Omega\bigr)$.
\end{lem}

\begin{proof}
The map descends to $\pi_0\Symp_{\rm fib}(X)$ using the invariance of Floer cohomology under Hamiltonian isotopy. Given $\phi_1$ and $\phi_2$ we can choose $\xi$ large enough that $T$, $\phi_1(T)$, $\phi_2(T)$ all lie inside~$X_\xi$; by arguments analogous to those above, we get three corresponding $(\Lambda_\bC^*)^2$-charts of point-like objects on $U_{\xi}$. The homomorphism property is then immediate. \end{proof}

We remind the reader that all these results are contingent on Assumption \ref{assump:inexact}.
But granted that,
by construction, the map of Lemma~\ref{lem:done} would split the homomorphism $\Bir_e\bigl(\bP^2,\Omega\bigr) \to \pi_0\Symp_{\rm fib}(X)$ given in \cite{Keating-Ward}.

\section{Further examples}\label{sec:sporadic}

\subsection{Tori}

For the standard symplectic structure on the torus $\bR^{2g}/\bZ^{2g}$, there is a linear action by symplectomorphisms, which is split by considering the action on first homology of the torus.

\subsection{Two-dimensional surfaces}

For a compact symplectic surface $\Sigma_g$, \cite{Auroux-Smith} shows that $\Gamma_g\to \Symp(\Sigma_g) / \mathrm{Ham}(\Sigma_g)$ splits by considering the composite of the monodromy homomorphism with the representation of $\Symp / \mathrm{Ham}$ on $\Auteq(\scrF(\Sigma_g))$, where we work over the Novikov field. That composite homomorphism splits; this is established by classifying spherical objects in the Fukaya category with non-trivial Chern character, and thereby reconstructing the `Schmutz graph' of the surface \cite{Schmutz}.

\subsection{(Log) Calabi--Yau surfaces}

Gromov \cite{Gromov-1985} showed that $\Symp_c(\bC^2, \omega_\std)$ is contractible. Later, Wendl \cite{Wendl} and Evans \cite{Evans} proved contractibility of $\Symp_c(\bC^* \times \bC^*)$ and $\Symp_c(\bC\times\bC^*)$ respectively.

Seidel \cite{Seidel_T*S2_unpublished, Seidel_lectures_4Dtwists} showed $\Symp_c T^\ast S^2$ is weakly homotopy equivalent to $\bZ$, and Evans \cite{Evans} proved that $\Symp_c\bigl(T^*\bR\bP^2\bigr) \simeq \bZ$, in both cases generated by the corresponding (spherical respectively projective) Dehn twist.

The symplectomorphism groups of log Calabi--Yau surfaces were studied systematically via homological mirror symmetry in \cite{Hacking-Keating-1,Hacking-Keating-2}. We summarise the relevant results.

Let $(Y, D)$ be a log Calabi--Yau surface with maximal boundary and split mixed Hodge structure, and let $M$ be its mirror. Let $U = Y \setminus D$. We have a mirror equivalence of categories $D(U) \cong D^b (\tw \scrW (M))$.
In this context, there is a natural subgroup of `compactly supported' auto-equivalences of $D(U)$, $\Auteq_c D(U)$ (see \cite[Definition~2.16]{Hacking-Keating-2}). Via the mirror equivalence, we have
\begin{equation*}
\pi_0 \Symp_c (M) \longrightarrow \Auteq_c D(U).
\end{equation*}

There exist conjectural descriptions of $\Auteq D(Y)$ and of $\Auteq_c D(U)$,
both in terms of generators and as fundamental groups of suitable moduli spaces of complex structures, see \cite[Conjecture~1.2]{Uehara}, \cite[Conjectures~1.5 and~1.7]{Hacking-Keating-2}. One can think of these as relative versions of the better-known conjecture of Bridgeland on autoequivalences of derived categories of $K3$ surfaces. In particular, we expect $\Auteq_c D(U) = \langle Q, \cT \rangle \rtimes \Aut (Y, D)$. Here $Q = \langle D_1, \dots, D_k \rangle^\perp \subset \Pic (Y)$, where $D_1, \dots, D_k$ are the irreducible components of $D$, is the subgroup of line bundles~$\cL$ on~$Y$ such that $\cL|_D = \cO_D$ (acting by tensor product); $\cT$ is the collection of spherical twists in objects of the form $i_\ast \cO_C (k)$, for $C \subset U$ a $(-2)$ curve and $k \in \bZ$; and $\Aut (Y, D )$ is the group of automorphisms of $Y$ fixing $D$ pointwise (acting by pushforward). By the Torelli theorem in \cite{GHK2}, $\Aut (Y, D )$ is identified with the quotient of the group $\Adm$ of admissible lattice automorphisms of $\Pic (Y)$ by the Weyl group~$W$.

The spherical objects $i_\ast \cO_C (k)$ are mirror to embedded Lagrangian spheres \cite[Proposition~5.2]{Hacking-Keating-2}, with spherical twists mirror to Dehn twists in the usual manner. Moreover, there is a~faithful group homomorphism \cite[Theorem~7.1]{Hacking-Keating-2}
\begin{equation} \label{eqn:logCYmonodromy}
\iota\colon \ Q \rtimes \Adm / W \hookrightarrow \pi_0 \Symp_c (M).
\end{equation}
The image symplectomorphisms are described in terms of the SYZ (almost-toric) Lagrangian fibration on $M$.
Under $\iota$, a line bundle $\cL$ acts as a `Lagrangian translation', with mirror autoequivalence $\otimes \cL \in \Auteq D(U)$; and an element $\phi \in \Adm / W$ acts as a `nodal slide recombination', with mirror autoequivalence $\phi_\ast \in \Auteq D(U)$.
In particular, combining back with Dehn twists, one has precise descriptions of the monodromy representation of the (conjectural) relevant moduli space; and, subject to the conjecture on $\Auteq_c D(U)$, we expect the representation~\eqref{eqn:logCYmonodromy} to split.

For projective Calabi--Yau surfaces, a general splitting result is expected to follow from Bridgeland's well-known conjecture on spaces of stability conditions on $K3$ surfaces.\footnote{In the abelian surface case, where there are no spherical objects and Bridgeland's conjecture is a theorem, there are few autoequivalences acting trivially on cohomology and the situation is close to that for product tori with splitting arising from lattice theory and the linear action of symplectomorphisms on $[\omega]^{\perp}$, cf. \cite{Fan-Lai}.}
Some Picard rank one $K3$ surfaces were studied via HMS in \cite{Sheridan-Smith:K3}. For a concrete example, one can take the `mirror quartic' $(X,\omega)$, which is a crepant resolution of the quotient of the Fermat quartic $\bigl\{\sum z_j^4 = 0\bigr\} \subset \bC\bP^3$ by an action of $(\bZ/4)^2$. The crepant resolution resolves six $A_3$-singularities, yielding 18 exceptional curves. Take a K\"ahler form which is `sufficiently irrational' in the sense that the areas of these 18 curves, and of a line, are all linearly independent over $\bQ$. The moduli space of complex structures is an orbifold curve with fundamental group
\[
\pi_1(\mathcal{M}_X) = \bZ \ast \bZ/4.
\]
The natural map $\pi_1(\mathcal{M}_X) \to \pi_0\Symp(X,\omega)$ indeed splits by \cite[Section 7.2]{Sheridan-Smith:K3}. Similar results are obtained for the crepant resolution of the Fermat sextic in $\bP(1,1,1,3)$, again with a sufficiently irrational K\"ahler form (meaning that the relevant complex structure moduli spaces is again only one-dimensional).

\begin{rmk}
 S.~Mu\~noz-Ech\'aniz proves an interesting splitting result at the level of \emph{abelianisations} of symplectic mapping class groups of $K3$ surfaces, cf.~\cite[Theorem B]{Munoz-Echaniz}. This is again based on the presence of infinite families of homologically distinct Lagrangian spheres.
\end{rmk}

\subsection{del Pezzo surfaces}

We consider rational surfaces which admit a monotone symplectic form, i.e.,~the surfaces $X_k$ which are monotone blow-ups of the projective plane at $0 \leq k \leq 8$ points (the del Pezzo of degree $9-k$) and the quadric surface $\bP^1 \times \bP^1$. The symplectomorphism groups in these cases have been extensively studied, and we summarise some of the results below. We write $\Symp_0$ for the Torelli group of symplectomorphisms acting trivially on homology. The algebro-geometric statements are taken from \cite{Dolgachev-book, Looijenga, Okada-Spotti-Sun}:
\begin{enumerate}\itemsep=0pt
 \item $\Symp_c\bigl(\bP^2, \omega_\FS\bigr) \simeq \bP {\rm U}(3)$.
 \item $\Symp_c\bigl(\bP^1\times \bP^1\bigr) \simeq ({\rm SO}(3)\times {\rm SO}(3)) \rtimes \bZ/2$, where $\bZ/2$ acts by swapping the factors.
 \item (degree 8) $\Symp \bigl(\Bl_1\bigl(\bP^2\bigr), \omega\bigl) \simeq {\rm U}(2)$.
 \item (degree 7) $\Symp\bigl(\Bl_2\bigl(\bP^2\bigr), \omega\bigr) \simeq T^2 \rtimes \bZ/2$, see \cite{Lalonde-Pinsonnault, Pinsonnault}.
 \item (degree 6) $\Symp_0\bigl(\Bl_3\bigl(\bP^2\bigr),\omega\bigr) \simeq T^2$ \cite{Evans}; the moduli space of complex structures is given by triples of distinct points in $\bP^2$ mod $\bP {\rm GL}_3(\bC)$.
 The holomorphic automorphism group of each surface is positive dimensional $(\bC^*)^2 \rtimes (S_3 \times S_2)$, with the finite part acting by permutations and quadratic Cremona transformations.
 \item (degree 5) $\Symp_0(\Bl_4(\bP^2)) \simeq \{pt\}$ \cite{Evans}. The moduli space is $\{pt\} / S_5$ (here $S_5 = \hat{A}_4$ also arises as a particular Artin group quotient in~\cite{Looijenga}). The monodromy action is split by taking the action on homology. See also
 \cite{Evans, Li-Li-Wu-1, Li-Li-Wu-2, Li-Wu}.
 \item (degree 4, intersections of two quadrics in $\bP^4$) We have $\Symp_0\bigl(\Bl_5\bigl(\bP^2\bigr)\bigr) \simeq \Diff^+\bigl(S^2,5\bigr) \simeq {\rm PBr}\bigl(S^2,5\bigr)/\langle \tau \rangle$, where $\tau$ is the full twist, which has order two. Seidel \cite{Seidel_lectures_4Dtwists} showed that the monodromy map for the universal homologically marked family defines a splitting
 \[
 {\rm PBr}\bigl(S^2,5\bigr)/\langle \tau \rangle \longrightarrow \pi_0\Symp_0\bigl(\Bl_5\bigl(\bP^2\bigr)\bigr) \longrightarrow {\rm PBr}\bigl(S^2,5\bigr)/\langle \tau \rangle.
 \]
 \item (degree 3, cubic surfaces) The fundamental group of the moduli space is a certain Artin group $\hat{E}_6$. As far as we know, this case has not been fully analysed, and splitting is unknown.

 \item (degree 2, quartics in $\bP^2$) Again, this case seems not to have been fully analysed.
\end{enumerate}
Separately, the reader may be interested in the related work \cite{Li-Li-Wu}, which inter alia applies to some non-monotone del Pezzo surfaces.

\subsection{A higher-dimensional non-splitting example?}

Let $X:= Q_0 \cap Q_1 \subset \bP^{2g+1}$ be the intersection of two smooth quadric hypersurfaces, with its canonical monotone symplectic form. Let $\mathcal{M}_X$ denote the moduli space of complex structures on $X$; these are all realised as $(2,2)$-complete intersections, by \cite[Theorem~1]{Ionescu}.

\begin{conj}
 For $g \geq 3$, the symplectic monodromy
$
 \pi_1(\mathcal{M}_X) \to \pi_0\Symp(X,\omega)$
 is not virtually split.
\end{conj}

 The two quadrics generate a pencil, whose singular members are generically parametrized by $2g+2$ points in $\bP^1$. This sets up an identification between the moduli space of complex structures on $X$ and a space of hyperelliptic curves, and $\mathcal{M}_X$ has fundamental group an extension $H^1(\Sigma_g;\bZ/2) \rtimes \Gamma_g^{\rm hyp}$ of the hyperelliptic mapping class group of a genus $g$ surface, see \cite{Reid,Smith}.
Smith in~\cite{Smith} constructs a natural map $ \pi_0\Symp(X) \to \Gamma_g$
 via autoequivalences of the Fukaya category. When $g=2$ the monodromy thus virtually splits, but for $g \geq 3$ the non-existence of a~map $\Gamma_g \to \Gamma_g^{\rm hyp}$ virtually splitting the inclusion suggests that the symplectic monodromy does \emph{not} split (more precisely, no splitting could factor through the autoequivalences of the Fukaya category). It would be interesting to investigate other Fano complete intersections.

 \begin{rmk}
 For degrees $d \geq 5$, the monodromy of the moduli space of plane curves is (at most) a proper subgroup of the mapping class group, corresponding to automorphisms which preserve an $r$-spin structure \cite{Salter}; conjecturally, the monodromy group is as large as this constraint permits, so whilst there is no splitting of the monodromy viewed as mapping to the classical mapping class group, there is a virtual splitting (one to a finite index subgroup). This example, as well as the case of the $A_k$-Milnor fibre, motivate considering virtual as well as actual splittings.\looseness=-1
\end{rmk}

\subsection[The C\^{}0-category]{The $\boldsymbol{C^0}$-category}

We end with a question. Following the work of Jannaud \cite{Jannaud}, in many cases one can show that a monodromy map $G \to \pi_0\Symp(X)$, whose injectivity is detected by computations of ranks of Floer cohomology groups, remains injective on further composition with the natural map
\[
\pi_0\Symp(X) \to \pi_0\overline{\Symp}(X).
\]
Here $\overline{\Symp}(X)$ refers to the group of symplectic homeomorphisms, that is, the $C^0$-closure of $\Symp(X)$ inside $\mathrm{Homeo}(X)$.

Such results are proved by exploiting the existence of a map $\overline{\Symp}(X) \to \mathcal{B}$ where $\mathcal{B}$ is a suitable space of (tuples of) barcodes with the bottleneck distance, and the homomorphism associates to a smooth symplectomorphism $\phi$ the barcode(s) associated to some particular sum of Floer complexes $\bigoplus_i {CF}(L_i,\phi(L_i))$ (where, for instance, the $\{L_i\}$ form a set of split-generators). The set of components $\pi_0(\mathcal{B})$ is not naturally a group.

\begin{qn}
 In the examples studied previously, does the monodromy homomorphism to the $C^0$-symplectic mapping class group ever split?
\end{qn}

\subsection*{Acknowledgements}
 A.K.~is partially supported by EPSRC Open Fellowship EP/W001780/1.
 I.S.~is partially supported by UKRI Frontier Research grant EP/X030660/1 (in lieu of an ERC Advanced Grant). M.W.~is partially supported by ERC Consolidator Grant 101001227. The authors are grateful to the anonymous referees for their numerous comments and suggested clarifications.

\subsection*{UKRI data access statement}
There is no dataset associated with this paper.


\pdfbookmark[1]{References}{ref}
\LastPageEnding

\end{document}